\documentstyle[11pt]{article}
\newcommand{\bbz}{\mbox{\boldmath $Z$}}
\newcommand{\bbn}{\mbox{\boldmath $N$}}
\newcommand{\bbp}{\mbox{\boldmath $P$}}

\newcommand{\bbq}{{\mathbf Q}}

\newcommand{\qed}{\hfill$\Box$}
\newcommand{\sn}{{\it span}}

\newcommand{\tr}{{\it Tr}}

\newcommand{\cont}{{\it cont}}
\newcommand{\inv}{{\it inv}}
\newcommand{\D}{{\it Des}}
\newcommand{\des}{{\it des}} 
\newcommand{\maj}{{\it maj}}
\newcommand{\DN}{{\it NDes}}
\newcommand{\N}{{\it Neg}}
\newcommand{\n}{{\it neg}}
\newcommand{\dn}{{\it ndes}}
\newcommand{\mn}{{\it nmaj}}
\newcommand{\fd}{{\it fdes}}  
\newcommand{\fm}{{\it fmaj}}
\newcommand{\la}{\lambda}
\newcommand{\hT}{{\hat{T}}}
\newcommand{\lala}{(\la^1,\la^2)}
\newcommand{\TT}{(T^1,T^2)}
\newcommand{\hTT}{(\hT^1,\hT^2)}

\newtheorem{thm}{Theorem}[section]
\newtheorem{pro}[thm]{Proposition}
\newtheorem{lem}[thm]{Lemma}
\newtheorem{cla}[thm]{Claim}

\newtheorem{cor}[thm]{Corollary}

\newtheorem{obs}[thm]{Observation}

\newenvironment{note}[1]{\par\addvspace{\medskipamount}\noindent
                         {\bf {#1}}\sl
                       }{\par\addvspace{\medskipamount}\rm}

\newcommand{\eqdef}{:=}

\begin{document}

\title{Descent Representations and Multivariate Statistics}
\bibliographystyle{acm}
\author{Ron M.\ Adin%
\thanks{Department of Mathematics and Computer Science, Bar-Ilan University,
Ramat-Gan 52900, Israel. Email: {\tt radin@math.biu.ac.il} }\ $^\S$ 
\and Francesco Brenti%
\thanks{Dipartimento di Matematica,
Universit\'{a} di Roma ``Tor Vergata'',
Via della Ricerca Scientifica,
00133 Roma, Italy. Email: {\tt brenti@mat.uniroma2.it} } $^\S$ 
\and Yuval Roichman%
\thanks{Department of Mathematics and Computer Science, Bar-Ilan University,
Ramat-Gan 52900, Israel. Email: {\tt yuvalr@math.biu.ac.il} } 
\thanks{Research of all authors was supported in part by 
the EC's IHRP programme, within the Research Training Network ``Algebraic
Combinatorics in Europe'', grant HPRN-CT-2001-00272,
by the Israel Science Foundation, founded by the Israel Academy of Sciences 
and Humanities, 
and by internal research grants from Bar-Ilan University.
}%
} 
\date{Submitted: October 13, 2002; Revised: June 12, 2003}

\maketitle

\begin {abstract}
Combinatorial identities on Weyl groups of types $A$ and $B$ are derived from 
special bases of the corresponding coinvariant algebras.
Using the Garsia-Stanton descent basis of the coinvariant algebra of type $A$ 
we give a new construction of the Solomon descent representations.
An extension of the descent basis to type $B$, using new multivariate statistics
on the group, yields a refinement of the descent representations.
These constructions are then applied to refine well-known decomposition rules 
of the coinvariant algebra and to generalize various identities.
\end{abstract}

\section{Introduction}

\subsection{Outline}

This paper studies the interplay between representations of classical Weyl groups 
of types $A$ and $B$ and combinatorial identities on these groups.
New combinatorial statistics on these groups are introduced, which lead to 
a new construction of representations. 
The Hilbert series which emerge give rise to multivariate identities 
generalizing known ones.

The set of elements in a Coxeter group having a fixed descent set carries 
a natural representation of the group, called a descent representation.
Descent representations of Weyl groups were first introduced by Solomon \cite{So} 
as alternating sums of permutation representations. This concept was extended to
arbitrary Coxeter groups, using a different construction, by Kazhdan and Lusztig
\cite{KL} \cite[\S 7.15]{Hum}. 
For Weyl groups of type $A$, these representations also appear
in the top homology of certain (Cohen-Macaulay) rank-selected posets \cite{St82}.
Another description (for type $A$) is by means of zig-zag diagrams \cite{Ge2, GR}. 

In this paper we give a new construction of descent representations for Weyl 
groups of type $A$, using the coinvariant algebra as a representation space. 
This viewpoint gives rise to a new extension for type $B$, which refines 
the one by Solomon.

The construction of a basis for the coinvariant algebra is important for 
many applications, and has been approached from different viewpoints.
A geometric approach identifies the coinvariant algebra with the cohomology ring
$H^*(G/B)$ of the flag variety.
This leads to the Schubert basis \cite{BGG, De}, and applies to any Weyl group.
This identification also appears in Springer's construction of irreducible 
representations \cite{Sp}; see also \cite{GP}. 
Barcelo \cite{Ba} found bases for the resulting quotients.
An algebraic approach, applying Young symmetrizers, was used by Ariki, 
Terasoma and Yamada \cite{TY, ATY} to produce a basis compatible with 
the decomposition into irreducible representations. 
This was extended to complex reflection groups in \cite{MY}.

A combinatorial approach, which produces a basis of monomials, was presented by
Garsia and Stanton in \cite{GS} (see also \cite{Gar, Stg}). 
They actually presented a basis for a finite dimensional quotient of the
Stanley-Reisner ring arising from a finite Weyl group. 
For type $A$, unlike other types, this quotient is isomorphic to the 
coinvariant algebra.
The Garsia-Stanton descent basis for type $A$ may be constructed from the 
coinvariant algebra via a straightening algorithm \cite{Al1}. 
Using a reformulation of this algorithm we give a natural 
construction of Solomon's descent representations as factors of the 
coinvariant algebra of type $A$.

An analogue of the descent basis for type $B$ is now given.
This analogue (again consisting of monomials) involves extended descent sets 
and new combinatorial statistics. 
An extension of the construction of descent representations, 
using the new basis for type $B$, 
gives rise to a family of descent representations, refining Solomon's. 
A decomposition of these descent representations into irreducibles,
refining theorems of Lusztig and Stanley  
\cite[Prop.~4.11]{St79} \cite[Theorem 8.8]{Re} (for type $A$)
and Stembridge \cite{Stem} (for type $B$), 
is carried out using a multivariate version of Stanley's formula for the 
principal specialization of Schur functions.  

This algebraic setting is then applied to obtain new multivariate 
combinatorial identities. 
Suitable Hilbert series are computed and compared to each other
and to generating functions of multivariate statistics.
The resulting identities present a far reaching generalization of 
bivariate identities from \cite{Ge1}, \cite{Gar79}, and \cite{ABR}.

\subsection{Main Results}

Let $W$ be a classical Weyl group of type $A$ or $B$, and let $I^W_n$ be 
the ideal of the polynomial ring
$P_n\eqdef \bbq[x_1,\dots,x_n]$ generated by $W$-invariant polynomials 
without a constant term.
The quotient $P_n/I^W_n$ is called the {\it coinvariant algebra} of $W$.
See Subsection 2.5 below.

For any partition $\la$ with (at most) $n$ parts, let 
$P_\la^{\underline{\triangleleft}}$ be the subspace of the polynomial 
ring $P_n=\bbq[x_1,\dots,x_n]$ spanned by all monomials whose exponent 
partition is dominated by $\la$, and
let $R_\la$ be a distinguished quotient of the image of 
$P_\la^{\underline{\triangleleft}}$ under the projection of 
$P_n$ onto the coinvariant algebra. 
For precise definitions see Subsections 3.5 and 5.3.
We will show that the homogeneous components of the coinvariant algebra 
decompose as direct sums of certain $R_\la$-s.
This will be done using an explicit construction of a basis for $R_\la$.
The construction of this basis involves new statistics on $S_n$ and $B_n$.

\subsubsection{New Statistics}

Let $\Sigma$ be a linearly ordered alphabet.
For any finite sequence $\sigma=(\sigma_1,\ldots,\sigma_n)$  
of letters in $\Sigma$ define
\[
\D(\sigma)\eqdef \{i\ |\ \sigma_i>\sigma_{i+1}\},
\]
the descent set of $\sigma$, and
\[
d_{i}(\sigma ) \eqdef |\{ j \in \D(\sigma ): \; j \geq i \} | ,
\]
the number of descents in $\sigma$ from position $i$ on.

\medskip

\noindent
If $\Sigma$ consists of integers, let
\[
\N(\sigma)\eqdef \{ i\ |\ \sigma_i<0\} ;
\]
\[
n_{i}(\sigma ) \eqdef
|\{ j \in \N(\sigma ): \; j \geq i \} | ;
\]
$$
 \varepsilon _{i} (\sigma ) \eqdef \left\{
\begin{array}{ll}  
1, & \mbox{if $\sigma_i< 0$,} \\
0, & \mbox{otherwise;}
\end{array} \right. 
$$
and
$$
f_i(\sigma):=2d_i(\sigma)+\varepsilon_i(\sigma).
$$
The statistics $f_i(\sigma)$ refine the flag-major index $\fm(\sigma)$,
which was introduced and studied in \cite{AR1, AR, ABR}.

For various properties of these statistics see Sections 3, 5, 
and 6 below.

\subsubsection{The Garsia-Stanton Descent Basis and its Extension}

To any $\pi\in S_n$ Garsia and Stanton \cite{GS} associated the monomial
$$
a_\pi\eqdef  \prod_{i\in \D(\pi)} (x_{\pi(1)}\cdots x_{\pi(i)}).
$$
It should be noted that in our notation 
$a_\pi=\prod_{i=1}^n x_{\pi(i)}^{d_i(\pi)}$.
Using Stanley-Reisner rings, Garsia and Stanton~\cite{GS} showed that 
the set $\{a_\pi+I_n\ |\ \pi\in S_n\}$ forms a basis for 
the coinvariant algebra of type $A$. 
This basis will be called the {\it descent basis}. The 
Garsia-Stanton approach is not applicable to the coinvariant 
algebras of other Weyl groups. In this paper we extend the 
 descent basis to the Weyl groups of type $B$.

To any $\sigma\in B_n$ we associate the monomial
$$
b_\sigma \eqdef
\prod_{i=1}^n x_{|\sigma(i)|}^{f_i(\sigma)}.
$$

\begin{thm}\label{main.0} [See Corollary~\ref{6.2'}]
The set 
$$
\{b_\sigma+I^B_n\ |\ \sigma\in B_n\}
$$
forms a basis for the coinvariant algebra of type $B$.
\end{thm}

\subsubsection{Descent Representations}

For a monomial $m$ in the polynomial ring $P_n=\bbq[x_1,\dots,x_n]$,
let the exponent partition $\la(m)$ be the partition obtained
by rearranging the exponents in a weakly decreasing order.
For any partition $\la$ with at most $n$ parts, let 
$P_\la^{\underline{\triangleleft}}$ be the subspace of  $P_n$ spanned by
all monomials whose exponent partition is dominated by $\la$ :
$$
P_\la^{\underline{\triangleleft}}:=\sn_{\bbq}\{m\ |\ \la(m)\ \underline{\triangleleft}\ \la\}.
$$
Similarly, define $P_\la^{\triangleleft}$ by strict dominance :
$$
P_\la^{\triangleleft}:=\sn_{\bbq}\{m\ |\ \la(m)\triangleleft\la\}.
$$
Consider now the canonical projection of $P_n$ onto the coinvariant algebra
$$
\psi : P_n \longrightarrow P_n/I_n.
$$
Define $R_\la$ to be a quotient of images under this map :
$$
R_\la \eqdef \psi(P_\la^{\underline{\triangleleft}})/\psi(P_\la^{\triangleleft}).
$$
Then $R_\la$ is an $S_n$-module. 

\smallskip

\noindent
For any subset $S\subseteq \{1,\dots,n\}$ define a partition
$$
\lambda_S \eqdef
(\lambda_1,\dots,\lambda_n)
$$
by
$$
\lambda_i \eqdef
|\,S\cap \{i,\dots,n\}\,|.
$$
Using a straightening algorithm for the descent basis it is shown that 
$R_\la\not=0$ if and only if $\la=\la_S$ for some $S\subseteq [n-1]$ (Corollary~\ref{R.la}), and that a basis for $R_{\la_S}$ may be indexed by 
the permutations with descent set equal to $S$ (Corollary~\ref{4.8}).
Let $R_k$ be the $k$-th homogeneous component
of the coinvariant algebra $P_n/I_n$.

\begin{thm}\label{basis.A} [See Theorem~\ref{a.main}]
For every $0\le k\le {n\choose 2}$, 
$$
R_k \cong \bigoplus_{S} R_{\la_S}
$$
as $S_n$-modules, where the sum is over all subsets
 $S\subseteq [n-1]$ such that $\sum_{i\in S} i =k$.
\end{thm}

\medskip

Let
$$
R^B_{\la}:=\psi^B(P^{\underline{\triangleleft}}_{\la})\ /\ 
\psi^B(P^{\triangleleft}_{\la}),
$$
where $\psi^B:P_n\longrightarrow P_n/I^B_n$ is the canonical map
from $P_n$ onto the coinvariant algebra of type $B$.
For subsets $S_1\subseteq [n-1]$
and $S_2\subseteq [n]$, let $\la_{S_1,S_2}$ be the vector
$$
\la_{S_1,S_2}:=2\la_{S_1}+{\mathbf 1}_{S_2},
$$
where $\la_{S_1}$ is as above and ${\mathbf 1}_{S_2}\in \{0,1\}^n$ is 
the characteristic vector of $S_2$. 
Again, $R^B_\la\not=0$ if and only if $\la=\la_{S_1,S_2}$ for some 
$S_1\subseteq [n-1], S_2\subseteq [n]$ (Corollary~\ref{R.la.b}). 
In this case, a basis for $R^B_{\la_{S_1,S_2}}$ may be indexed by 
the signed permutations $\sigma\in B_n$ with $\D(\sigma)=S_1$ 
and $\N(\sigma)=S_2$ (Corollary~\ref{4.8.b}).

Let $R^B_k$ for the $k$-th homogeneous component
of the coinvariant algebra of type $B$.
The following theorem is a $B$-analogue of Theorem \ref{basis.A}.

\begin{thm}\label{basis.B} [See Theorem~\ref{b.4}]
For every $0\le k\le n^2$,
$$
R^B_k \cong \bigoplus_{S_1,S_2} R^B_{\la_{S_1,S_2}}
$$
as $B_n$-modules, where the sum is over all subsets
$S_1\subseteq [n-1]$ and $S_2\subseteq [n]$ such that
$\lambda_{S_1,S_2}$ is a partition and
$$
2\cdot \sum_{i\in S_1} i +|S_2|=k.
$$
\end{thm}

\subsubsection{Decomposition into Irreducibles}

\begin{thm}\label{main.3} [See Theorem~\ref{t.rep.decomp}]
For any subset $S\subseteq [n-1]$ and partition $\mu\vdash n$,
the multiplicity in $R_{\la_S}$ of the irreducible $S_n$-representation
corresponding to $\mu$ is
$$
m_{S,\mu} \eqdef
|\,\{\,T\in SYT(\mu)\ |\ \D(T)=S\,\}\,|, 
$$
the number of standard Young tableaux of shape $\mu$ and descent set $S$.
\end{thm}

\noindent
This theorem refines the well known decomposition
(due, independently, to Lusztig and Stanley)
of each homogeneous component of the coinvariant algebra into
irreducibles. See Subsection 2.5.

\medskip

For type $B$ we have

\begin{thm}\label{main.4}
For any pair of subsets $S_1\subseteq [n-1]$, $S_2\subseteq [n]$,
and a bipartition $(\mu^1,\mu^2)$ of $n$, the multiplicity of
the irreducible $B_n$-representation corresponding to $(\mu^1,\mu^2)$
in $R^B_{\la_{S_1,S_2}}$ is
$$
m_{S_1,S_2,\mu^1,\mu^2} \eqdef
|\,\{\,T\in SYT(\mu^1,\mu^2)\ | \; \D(T)=S_1 \hbox{ and } \N(T)=S_2\,\}\,|, 
$$
the number of pairs of standard Young tableaux of shapes $\mu^1$ and $\mu^2$ 
with descent set $S_1$ and sets of entries $[n]\setminus S_2$ and $S_2$, 
respectively.  
\end{thm}

\noindent
For definitions and more details see Subsection 5.4 and 
Theorem~\ref{t.rep.decomp.B} below.

\medskip

The proofs apply multivariate extensions of Stanley's formula
for the principal specialization of a Schur function.
See Lemmas \ref{t.tau} and \ref{t.tau.B} below.

\subsubsection{Combinatorial Identities}

For any partition $\la = (\la_1,\ldots,\la_n)$ with at most $n$ parts define
$$
m_j(\la) \eqdef
|\ \{1\le i\le n\ |\ \la_i=j\}\ |\qquad (\forall j\ge 0).
$$

By considering Hilbert series of the polynomial ring with respect to rearranged multi-degree and applying the Straightening Lemma for the coinvariant algebra 
of type $A$ we obtain

\begin{thm}\label{comb.A} [See Theorem~\ref{7.2}]
For any positive integer $n$
\[ 
\sum _{\ell(\lambda)\le n} 
{n\choose m_0(\lambda),m_1(\lambda),\dots} \prod_{i=1}^n q_i^{\lambda_i} = 
\frac{{
\sum _{\pi \in S_{n}} \prod_{i=1}^n q_i^{d_i(\pi)}
}}
{{
 \prod_{i=1}^{n}}(1-q_1\cdots q_i)}
\]
in ${\bbz}[[q_1, \ldots , q_n]]$, where the sum on the left-hand side is taken over all
partitions with at most $n$ parts. 
\end{thm}

\noindent
This theorem generalizes Gessel's theorem for the bivariate distribution
of descent number and major index \cite{Ge1}.

\medskip

The main combinatorial result for type $B$ asserts :

\begin{thm}\label{main.comb} [See Theorem~\ref{7.1}]
For any positive integer $n$ 
\[ \sum _{\sigma \in B_{n}} \prod_{i=1}^n 
q_i^{d_i(\sigma)+n_i(\sigma^{-1})} = \sum _{\sigma \in B_{n}}\prod_{i=1}^n 
q_i^{2d_i(\sigma)+\varepsilon_i(\sigma)}.
 \]
\end{thm}

For further identities see Section 6.
In particular, it is shown that central results from \cite{ABR} follow from 
Theorem~\ref{main.comb}.

\section{Preliminaries}

\subsection{Notations}

Let $\bbp \eqdef \{ 1,2,3, \ldots \}$, $\bbn \eqdef\ \bbp \cup\ \{ 0 \}$, 
$\bbz$ be the ring of integers, and $\bbq$ be the field of rational numbers;
for $a \in \bbn$ let $[a] \eqdef \{ 1,2, \ldots , a \} $ (where $[0] \eqdef \emptyset $). 
Given $n, m \in \bbz$, $n \leq m$, let $[n,m] \eqdef \{ n,n+1, \ldots ,m \}$. 
For $S \subset \bbn$ write $S= \{ a_{1}, \ldots , a_{r} \} _{<}$ to mean that 
$S= \{ a_{1}, \ldots , a_{r} \}$ and $a_{1} < \ldots < a_{r} $.
The cardinality of a set $A$ will be denoted by  
$|A|$. More generally, given
 a multiset $M =  \{ 1^{a_{1}},2^{a_{2}},\ldots , r^{a_{r}} \}$
denote by $|M|$  its cardinality, so $|M|=\sum_{i=1}^{r}a_{i}$. 

Given a variable $q$ and a commutative ring $R$, denote by $R[q]$ (respectively,  $R[[q]]$) 
the ring of polynomials (respectively, formal power series) in $q$ with coefficients in $R$. For $i \in \bbn$
let, as customary, $[i]_{q} \eqdef 1+q+q^{2}+\ldots + q^{i-1}$ (so $[0]_{q}=0$).
Given a vector $v=(v_1,\dots,v_n)\in \bbp^n$ and a sequence of variables 
$x_1,\dots,x_n$ denote by $\bar x^v$ the monomial $\prod_{i=1}^n x_i^{v_i}\in R[x_1,\dots,x_n]$.

\subsection{Sequences and Permutations}

Let $\Sigma$ be a linearly ordered alphabet.
Given a sequence $\sigma =(\sigma_{1}, \ldots ,\sigma_{n}) \in \Sigma^{n}$ 
we say that a pair $(i,j) \in [n] \times [n]$ is an {\em inversion} of $\sigma $ if $i<j$ and $\sigma_{i}>\sigma_{j}$.
We say that $i \in [n-1]$ is a {\em descent} of $\sigma $ if $\sigma_{i}>\sigma_{i+1}$ ;
$$
\D(\sigma) \eqdef \{1\le i\le n-1\ | \; \sigma_{i} > \sigma_{i+1} \}
$$
is the {\em descent set} of $\sigma$.
Denote by $\inv(\sigma )$ (respectively, $\des(\sigma )$) the number of inversions (respectively, descents) of
$\sigma $. 
We also let
\[ maj (\sigma ) \eqdef \sum _{i\in\D(\sigma)}  i \]
and call it the {\em major index} of $\sigma$.

Given a set $T$ let $S(T)$ be the set of all bijections $\pi : T \rightarrow T$, and $S_{n} \eqdef S([n])$.
For $\pi \in S_{n}$ write $\pi = \pi_{1} \ldots \pi_{n}$ to mean that $\pi (i) =\pi_{i}$, for $i=1, \ldots ,n$. 
Any $\pi \in S_{n}$ may also be written in {\em disjoint cycle form} (see, e.g., \cite[p.17]{StEC1}), usually omitting
the 1-cycles of $\pi$. For example, $\pi =365492187$ may also be written as $\pi =(9,7,1,3,5)(2,6)$. 
Given $\pi , \tau \in S_{n}$ let $\pi \tau \eqdef \pi \circ \tau $ (composition of functions) 
so that, for example, $(1,2)(2,3)=(1,2,3)$.

Denote by $B_{n}$ the group of all bijections $\sigma$ of the set $[-n,n]
\setminus \{ 0 \}$ onto itself such that 
\[ \sigma (-a)=-\sigma (a) \]
for all $a \in [-n,n] \setminus \{ 0 \}$, with composition as the group operation.
This group is usually known as the group of ``signed permutations''
on $[n]$, or as the {\em hyperoctahedral group} of rank $n$. We identify $S_{n}$
as a subgroup of $B_{n}$, and $B_{n}$ as a subgroup of $S_{2n}$, in the
natural ways.

For $\sigma \in B_{n}$ write $\sigma = [a_{1}, \ldots ,a_{n}]$ to mean that $\sigma (i)=a_{i}$ for $i=1, \ldots , n$, 
and (using the natural linear order on $[-n,n]\setminus \{0\}$) let
$$
\begin{array}{ll}
\inv(\sigma) \eqdef \inv (a_{1}, \ldots , a_{n} ) , &
  \des(\sigma) \eqdef \des(a_{1}, \ldots , a_{n}) ,\\
\N(\sigma) \eqdef \{ i \in [n]: \; a_{i}<0 \} ,&
  \n(\sigma) \eqdef |\N(\sigma )| , \\
\maj(\sigma) \eqdef \maj (a_{1}, \ldots , a_{n} ) , &
  \fm(\sigma) \eqdef 2\cdot \maj(\sigma)+ \n(\sigma).
\end{array}
$$
The statistic $\fm$ was introduced in \cite{AR1, AR} and further studied in 
\cite{ABR}.

\subsection{Partitions and Tableaux}\label{s.part_tab}

Let $n$ be a nonnegative integer. A {\it partition} of $n$
is an infinite sequence of nonnegative integers with finitely many nonzero terms
$\lambda=(\lambda_1,\lambda_2,\ldots)$,
where $\lambda_1\ge\lambda_2\ge\ldots$ and  $\sum_{i=1}^{\infty} \lambda_i =n$. 
The sum $\sum \la_i=n$ is called the {\it size} of $\lambda$, denoted $|\la|$; write also $\lambda\vdash n$. 
The number of parts of $\la$, $\ell(\la)$, is the maximal $j$ for which $\la_j>0$.
The unique partition of $n=0$ is the {\it empty partition} $\emptyset=(0,0,\dots)$, 
which has length $\ell(\emptyset):=0$.
For a partition  $\lambda=(\lambda_1,\ldots,\lambda_k,\ldots)$  define the {\it conjugate partition}
$\lambda'=(\lambda'_1,\dots,\lambda'_i,\ldots)$ by letting $\lambda'_i$ be the 
number of parts of $\lambda$ that are $\ge i$ $(\forall i\ge 1)$. 


The {\it dominance} partial order on partitions is defined as follows :        
For any two partitions $\mu$ and $\la$ of the same integer,
$\mu$ {\it dominates} $\la$ (denoted $\mu\ \underline{\triangleright}\ \la$)
if and only if $\sum_{j=1}^i \mu_j \ge  \sum_{j=1}^i \la_j$ for all $i$ (and, by assumption, 
$\sum_{j=1}^\infty \mu_j = \sum_{j=1}^\infty \la_j$).


The subset $\{(i,j)\ |\ i,j\in {\mathbf P} ,\ j\le\la_i\}$ of ${\mathbf P}^2$ is called the {\it Young diagram}
of {\it shape} $\la$.  $(i,j)$ is the {\it cell} in row $i$ and column $j$.
The diagram of the conjugate shape $\la'$ may be obtained from the diagram of shape $\la$ 
by interchanging rows and columns.


A {\it Young tableau} of shape $\lambda$ is obtained by inserting the integers 
$1,\ldots,n$ (where $n=|\la|$) as {\it entries} in the cells of the Young 
diagram of shape $\lambda$, allowing no repetitions.
A {\it standard Young tableau} of shape $\lambda$ is 
a Young tableau whose entries increase along rows and columns.


We shall draw Young tableaux as in the following example.

\smallskip

\noindent{\bf Example 1.}
\[
\begin{array}{ccccc}
1         & 3         & 4                 & 6     & 9     \\
2         & 7         & 8                 & 11        \\
5         & 10          
\end{array}
\]

\smallskip

A {\it descent} in a standard Young tableau $T$ is an entry
$i$ such that $i+1$ is strictly south (and hence weakly west) of $i$.
Denote the set of all descents in $T$ by $\D(T)$.
The {\it descent number} and the {\it major index} (for tableaux) 
are defined as follows :
$$
\des(T)\eqdef\sum_{i\in \D(T)} 1\hbox{ };\qquad
\maj(T)\eqdef\sum_{i\in \D(T)} i.
$$

\noindent{\bf Example 1. (Cont.)} 
Let $T$ be the standard Young tableau drawn in Example 1.
Then $\D(T)=\{1,4,6,9\}$, $\des(T)=4$, and $\maj(T)=1+4+6+9=20$.

\medskip

A {\it semistandard Young tableau} of shape $\lambda$
is obtained by inserting positive integers as entries in the
cells of the Young diagram of shape $\lambda$,  
so that the entries weakly increase along rows and 
strictly increase down columns.
A {\it reverse semistandard Young tableau} is obtained by inserting 
positive integers into the diagram so that the entries weakly decrease 
along rows and strictly decrease down columns.

\smallskip

A {\it bipartition} of $n$ is a pair $(\la^1,\la^2)$ of partitions of total size $|\la^1|+|\la^2|=n$. 
A {\it (skew) diagram} of shape $(\la^1,\la^2)$ is the disjoint union of 
a diagram of shape $\la^1$ and a diagram of shape $\la^2$,
where the second diagram lies southwest of the first.
A standard Young tableau $T=\TT$ of shape $(\la^1,\la^2)$ is obtained by 
inserting the integers $1,2,\ldots,n$
as entries in the cells, such that the entries increase along rows and columns.  
The descent set $\D(T)$, the descent number $\des(T)$,
and the major index $\maj(T)$
 of $T$ are defined as above.
The {\it negative set}, $\N(T)$, of such a tableau $T$ is the 
set of entries in the cells of $\la^2$.  
Define $\n(T):=|\N(T)|$ and
$$
\fm(T):=2\cdot \maj(T)+ \n(T).
$$

\medskip

\noindent{\bf Example 2.}  Let $T$ be

\[
\begin{array}{ccccc}
  &   & 2 & 5 & 6 \\
  &   & 3 &   &   \\
1 & 7 &   &   &   \\
4 & 8 &   &   &
\end{array}
\]

\noindent
We shall also consider $T$ as a pair of tableaux
$$
\TT=\left(
\begin{array}{ccc}
2 & 5 & 6 \\
3 &   & 
\end{array}
\;,\;
\begin{array}{cc}
1 & 7 \\
4 & 8 
\end{array}
\right).
$$

\noindent
$T$ is a standard Young tableau of shape $((3,1),(2,2))$;
$\D(T)= \{2,3,6,7\}$, 
$\maj(T)=18$, 
$Neg(T)=\{1,4,7,8\}$, $\n(T)=4$, and $\fm(T)=40$.

\medskip

Denote by $SYT(\la)$  the set of all standard Young tableaux of shape $\la$.
Similarly, $SSYT(\la)$ and $RSSYT(\la)$ denote the sets of all semistandard 
and reverse semistandard Young tableaux of shape $\la$, respectively.
$SYT(\la^1,\la^2)$ denotes the set of all standard Young tableaux of shape 
$(\la^1,\la^2)$.

\medskip

Let $\hat T$ be a (reverse) semistandard Young tableau. 
The {\it content vector} $\cont(\hat T)=(m_1,m_2,\dots)$ is defined by
$$
m_i\eqdef 
|\{\hbox{cells in $\hT$ with entry $i$}\}| \qquad (\forall i\ge 1).
$$


\subsection{Symmetric Functions}

Let $\Lambda$ ($\Lambda_n$)  be the ring of symmetric functions 
in infinitely many (resp. $n$)
variables.
In this paper we consider three bases for the vector space $\Lambda$:

\noindent
The {\it elementary symmetric functions} are defined as follows. 
For a nonempty partition $\la=(\la_1,\la_2,\dots,\la_t)$
$$
e_\la(\bar x)\eqdef \prod_{i=1}^t e_{\la_i}(\bar x),
$$
where $\bar x=(x_1,x_2,\dots)$ and
for any $k\in \bbp$
$$
e_k(\bar x)\eqdef \sum_{i_1<i_2<\cdots<i_k} x_{i_1}\cdots x_{i_k}.
$$
The {\it power sum symmetric functions} are defined as follows. 
For a nonempty partition $\la=(\la_1,\la_2,\dots,\la_t)$
$$
p_\la(\bar x)\eqdef \prod_{i=1}^t p_{\la_i}(\bar x),
$$
where for any $k\in \bbp$
$$
p_k(\bar x)\eqdef \sum_{i=1}^\infty x_i^k.
$$
The {\it Schur functions} are defined in a different manner.
For a nonempty partition $\la$ define
$$
s_\lambda(\bar x)
\eqdef\sum_{T\in SSYT(\lambda)} \bar x^{\cont(T)}.
$$
For the empty partition $\la=\emptyset$, define $e_\emptyset(\bar x)
=p_\emptyset(\bar x)=s_\emptyset(\bar x):=1$.

\begin{cla}\label{t.schur} (See \cite[Prop.\ 7.10.4]{StEC2})
$$
s_\la(\bar x)=\sum_{T\in RSSYT(\la)} \bar{x}^{\cont(T)}.
$$
\end{cla}

\noindent
Note that corresponding spanning sets for $\Lambda_n$ are obtained by 
substituting $x_{n+1}=x_{n+2}=\cdots=0$.

\medskip

The conjugacy classes of the symmetric group $S_n$
are described by their cycle type; thus, by partitions of $n$.
The irreducible representations of $S_n$ are also indexed
by these partitions (cf.~\cite{Sa}).

Let $\lambda$ and $\mu$ be partitions of $n$.
Denote by $\chi^\lambda_\mu$ the value, at a conjugacy class of cycle type
$\mu$, of the character of the irreducible $S_n$-representation corresponding 
to $\lambda$. These character values are entries in the transition matrix 
between two of the above bases. 
This fact was discovered and applied by Frobenius as an efficient tool for 
calculating characters (cf.~\cite[p. 401]{StEC2}).

\begin{note}
\noindent
{\bf Frobenius Formula.} \cite[Corollary 7.17.4]{StEC2}
For any partition $\mu$ of $n$
$$
p_{\mu}(\bar x)=\sum_{\lambda\vdash n} \chi^{\lambda}_{\mu} s_\lambda(\bar x),
$$
where the sum runs through all partitions $\lambda$ of $n$.
\end{note}

\medskip

Recall that the Weyl group $B_n$
may be described as the group of all {\it signed permutations} of order $n$
(Subsection 2.2).
The conjugacy classes of $B_n$ are described by the signed cycle type 
(i.e., each cycle is described by its length and the product of the signs of 
its entries). Thus, the conjugacy classes are described by ordered pairs of 
partitions whose total size is $n$: the first partition consist of the lengths 
of the positive cycles, while the second consists of the lengths of the 
negative cycles.
The irreducible representations of $B_n$ are also indexed by these pairs 
of partitions, cf. \cite[\S I, Appendix B]{Md}.

\smallskip

For two partitions $\mu^1$ and $\mu^2$, and two
infinite sets of independent variables $\bar x=x_1,x_2,\dots$ 
and $\bar y=y_1,y_2,\dots$, define
$$
p_{\mu^1,\mu^2}(\bar x,\bar y)\eqdef 
\prod_i (p_{\mu^1_i}(\bar x)+p_{\mu^1_i}(\bar y))\cdot
\prod_j (p_{\mu^2_j}(\bar x)-p_{\mu^2_j}(\bar y)).
$$
Let $\chi^{\lambda^1,\lambda^2}_{\mu^1,\mu^2}$ be the 
irreducible $B_n$-character indexed by $(\lambda^1,\lambda^2)$ 
evaluated at a conjugacy class
of cycle type $(\mu^1,\mu^2)$.
Then

\begin{note}
\noindent
{\bf Frobenius Formula for $B_n$.} \cite[p. 178]{Md}
With the above notations, for any bipartition $(\mu^1,\mu^2)$ of $n$
$$
p_{\mu^1,\mu^2}(\bar x,\bar y)=
\sum_{\la^1,\la^2} \chi^{\la^1,\la^2}_{\mu^1,\mu^2} 
                   s_{\la^1}(\bar x) s_{\la^2} (\bar y),
$$
where the sum runs through all ordered pairs of partitions $(\la^1,\la^2)$
of total size $|\la^1|+|\la^2|=n$.
\end{note}

\subsection{The Coinvariant Algebra}

The groups $S_n$ and $B_n$ have natural actions on the ring of polynomials $P_n$ 
(cf. \cite[\S 3.1]{Hum}).
$S_n$ acts by permuting the variables, and $B_n$ acts by permuting the variables
and multiplying by $\pm 1$.
The ring of $S_n$-invariant polynomials is $\Lambda_n$, the ring of symmetric 
functions in $x_1,\dots, x_n$.
Similarly, the ring of $B_n$-invariant polynomials is $\Lambda^B_n$, the ring 
of symmetric functions in $x_1^2,\dots, x_n^2$.
Let $I_n$, $I^B_n$ be the ideals of $P_n$ generated by the elements
of $\Lambda_n$, $\Lambda^B_n$ (respectively) without constant term.
The quotient $P_n/I_n$ ($P_n/I^B_n$) is called the {\it coinvariant algebra} 
of $S_n$ ($B_n$).
Each group acts naturally on its coinvariant algebra. 
The resulting representation is isomorphic to the regular representation.
See, e.g., \cite[\S 3.6]{Hum} and \cite[\S II.3]{Hi}.

Let $R_k$ ($0\le k \le {n\choose 2}$) be the $k$-th homogeneous component of the coinvariant algebra 
of $S_n$: $P_n/I_n=\oplus_k R_k$. 
Each $R_k$ is an $S_n$-module.
The following theorem is apparently due to G.\ Lusztig (unpublished)
and, independently, to R.\ Stanley~\cite[Prop.~4.11]{St79}. 
It was also proved by Kraskiewicz and Weyman~\cite{KW}; see~\cite[p. 215]{Re}.

\begin{note}
\noindent{\bf Lusztig-Stanley Theorem.} 
\cite[Prop.~4.11]{St79} \cite[Theorem 8.8]{Re} 
For any $0\le k \le {n\choose 2}$ and $\mu\vdash n$,
the multiplicity in $R_k$ of the irreducible $S_n$-representation
corresponding to $\mu$ is
$$
m_{k,\mu}=
|\ \{\ T\in SYT(\mu)\ |\ \maj(T)=k\ \}\ | .
$$
\end{note}


\noindent
The following $B$-analogue (in different terminology) was proved in \cite{Stem}.
Here $R^B_k$ is 
the $k$-th homogeneous component of the coinvariant algebra of $B_n$.

\begin{note}
\noindent{\bf Stembridge's Theorem.}  For any 
$0\le k \le n^2$ and bipartition $(\mu^1,\mu^2)$ of $n$,
the multiplicity in $R^B_k$ of the irreducible $B_n$-representation
corresponding to $(\mu^1,\mu^2)$ is
$$
m_{k,\mu^1,\mu^2}=
|\ \{\ T\in SYT(\mu^1,\mu^2)\ |\ \fm(T)=k\ \}\ | .
$$
\end{note}

\subsection{The Garsia-Stanton Descent Basis (Type $A$)}

For any $\pi\in S_n$ define the monomial 
$$
a_\pi\eqdef \prod_{j\in\D(\pi)} (x_{\pi(1)}\cdots x_{\pi(j)}).
$$
Using Stanley-Reisner rings Garsia and Stanton showed that the set 
$\{a_\pi+I_n\ |\ \pi\in S_n\}$ forms a basis for the 
coinvariant algebra of type $A$ \cite{GS}. This basis 
will be called the {\it descent basis}. 
Unfortunately, this 
approach does not give a basis for the coinvariant algebras 
of other Weyl groups. In this paper we shall find an analogue of the  
descent basis for $B_n$.
Our approach will allow us to obtain refinements of the 
Lusztig-Stanley and Stembridge theorems.

\section{Construction of Descent Representations}

In this section we introduce several combinatorial concepts
(Subsections 3.1 and 3.2).
Then we present a straightening algorithm for the expansion of 
an arbitrary monomial in $P_n$ in terms of descent basis elements,
with coefficients from $\Lambda_n$
(Subsections 3.3 and 3.4). 
This algorithm is essentially equivalent to the one presented by 
Allen~\cite{Al1}, but our formulation leads naturally to 
descent representations and may be extended to type $B$ 
(Subsections 3.5 and 5.2).

\subsection{Descent Statistics}

In this subsection we introduce a family of descent statistics for
sequences of letters from a linearly ordered alphabet.

Let $\Sigma$ be a linearly ordered alphabet.
For any sequence $\sigma$  of $n$ letters from $\Sigma$ define
\[
d_{i}(\sigma ) \eqdef |\{ j \in \D(\sigma ): \; j \geq i \} | \qquad(1\le i\le n),
\]
the number of descents in $\sigma$ from position $i$ on.
Clearly,
$$
d_1(\sigma)=\des(\sigma)\leqno(3.1)
$$
and
$$
\sum_{i=1}^n d_i(\sigma) = \maj(\sigma). \leqno(3.2)
$$
Note that $d_n(\sigma)=0$.

Also,
for any sequence $\sigma$ from $\Sigma$
$$
d_i(\sigma)\ge d_{i+1}(\sigma), \qquad 1\le i < n ; \leqno(3.3)
$$
$$
d_i(\sigma)=d_{i+1}(\sigma)\Longrightarrow \sigma(i)<\sigma(i+1) ; \leqno(3.4)
$$
$$
d_i(\sigma)>d_{i+1}(\sigma) \Longrightarrow 
\sigma(i)>\sigma(i+1) \hbox{ and } d_i(\sigma)=d_{i+1}(\sigma)+1 . \leqno(3.5)
$$

It follows from (3.2) and (3.3) that
the sequence $(d_1(\sigma),d_2(\sigma),\dots,d_n(\sigma))$
is a partition of $\maj(\sigma)$. 
It follows from (3.4) and (3.5) that this sequence uniquely determines the 
descent set of $\sigma$. 
More explicitly, 
the  partition conjugate to
$(d_1(\sigma),\dots,d_n(\sigma))$ is obtained
by writing the elements of $\D(\sigma)$
in decreasing order.




Finally, for a permutation $\pi \in S_n$, let $a_\pi$ be 
the corresponding descent basis element. Then, clearly,
$$
a_\pi=\prod_{i=1}^n x_{\pi(i)}^{d_i(\pi)}.
$$

\subsection{\bf The Index Permutation and the Exponent Partition}
\label{s.index_exp}

The {\it index permutation} of a monomial $m=\prod_{i=1}^n x_i^{p_i}\in P_n$ 
is the unique permutation $\pi=\pi(m) \in S_n$ such that
$$
p_{\pi(i)}\ge p_{\pi(i+1)} \qquad (1\le i < n)  \leqno(1)
$$
and
$$
p_{\pi(i)}=p_{\pi(i+1)} \Longrightarrow \pi(i)<\pi(i+1). \leqno(2)
$$
In other words, $\pi$ reorders the variables $x_i$ by (weakly) decreasing exponents, where the 
variables with a given exponent are ordered by increasing indices.

\smallskip

Let $m=\prod_{i=1}^n x_i^{p_i}$ be a monomial in $P_n$,
$\pi=\pi(m)$ its index permutation, and $a_\pi$ the corresponding 
descent basis element. 
Then $m=\prod_{i=1}^n x_{\pi(i)}^{p_{\pi(i)}}$ and 
$a_\pi=\prod_{i=1}^n x_{\pi(i)}^{d_i(\pi)}$.

\begin{cla}\label{t.exp}
The sequence $(p_{\pi(i)}-d_i(\pi))_{i=1}^n$ of exponents in $m/ a_\pi$ 
consists of nonnegative integers, and is weakly decreasing:
$$
p_{\pi(i)}-d_i(\pi)\ge p_{\pi(i+1)}-d_{i+1}(\pi) \qquad (1\le i< n) .
$$
\end{cla}

\noindent{\bf Proof.} If $\pi(i)<\pi(i+1)$ then
$d_i(\pi)=d_{i+1}(\pi)$ and the claim follows from $p_{\pi(i)}\ge p_{\pi(i+1)}$. 
If $\pi(i)>\pi(i+1)$ then
$d_i(\pi)= d_{i+1}(\pi)+1$ and $p_{\pi(i)}>p_{\pi(i+1)}$,
so that $p_{\pi(i)}-p_{\pi(i+1)}\ge
1=d_i(\pi)-d_{i+1}(\pi)$.
This proves monotonicity.
Nonnegativity now follows from $d_n(\pi)=0$, implying 
$p_{\pi(n)}-d_n(\pi) \ge 0$.

\qed

\smallskip

For a monomial $m=\prod_{i=1}^n x_i^{p_i}$ with index permutation
$\pi\in S_n$, let 
$$
\lambda(m):=(p_{\pi(1)},p_{\pi(2)},\dots,p_{\pi(n)})
$$
be its {\it exponent partition}. Note that $\lambda(m)$ is a partition of 
the total degree of $m$.

 Define the {\it complementary partition}
$\mu(m)$ of a monomial $m$ to be the partition conjugate to the partition 
$(p_{\pi(i)}-d_i(\pi))_{i=1}^n$. Namely,
$$
\mu_j \eqdef |\{i\ |\ p_{\pi(i)}-d_i(\pi)\ge j\}| \qquad (\forall j\ge 1).
$$

\smallskip

\noindent{\bf Example.} Let $m=x_1^2x_2^4x_3^2x_5x_6^3$ and $n=7$.

\noindent
Then $m=x_2^4x_6^3x_1^2x_3^2x_5^1x_4^0x_7^0$, $\lambda(m)=(4,3,2,2,1,0,0)$,
$\pi=2613547\in S_7$,
$(d_1(\pi),\dots,d_7(\pi))=(2,2,1,1,1,0,0)$, 
and $\mu(m)= (4,1)$.



\smallskip

\subsection{A Partial Order for Monomials}

We shall now define a partial order on the monomials in $P_n$.

\noindent
Comparable monomials will always have the same total degree.
Fix such a total degree $p$; we may assume $p\ge 1$.

\smallskip

\noindent{\bf Definition.} For monomials $m_1,m_2$ of the same total degree
$p$ , $m_1\prec m_2$ if one of the following holds :
\begin{itemize}
\item[(1)] $\lambda(m_1)\ \triangleleft\ \lambda(m_2)$ (strictly smaller in dominance order); or
\item[(2)] $\lambda(m_1) = \lambda(m_2)$ and $\inv (\pi(m_1)) > \inv (\pi(m_2))$.
\end{itemize}

The partial order ``$\prec$" may be replaced, in this paper, by 
any linear extension of it;
for example, by the linear extension obtained by replacing the dominance order
by lexicographic order in (1).
(This is common  in the definition of Gr\"obner bases, which partially motivated our definition.)
  
\subsection{Straightening}

\begin{lem}\label{t.mk}
Let $m\in P_n$ be a monomial and
let $1\le k\le n$.
Let $S^{(k)}$ be the set of all monomials which appear (with coefficient 1)
in the expansion of the polynomial $m\cdot e_k$,
the product of $m$ and the $k$-th elementary symmetric function.
Let 
$$
m^{(k)}\eqdef m\cdot x_{\pi(1)}\cdot x_{\pi(2)}\cdots x_{\pi(k)},
$$
where $\pi=\pi(m)$ is the index permutation of $m$.
Then :
\begin{itemize}
\item[(1)] $m^{(k)}\in S^{(k)}$.
\item[(2)] $m'\in S^{(k)} \hbox{ and } m'\not=m^{(k)} \Longrightarrow m'\prec m^{(k)}$ .
\end{itemize}
\end{lem}

\noindent{\bf Proof.}
(1) is obvious. For (2), let $m'\in S^{(k)}$. Then
$m'=m\cdot x_{i_1}x_{i_2}\cdots x_{i_k}$, where $1\le i_1< i_2< \cdots < i_k \le n$.
Now, let $\la_i(m)$ be the $i$-th part of the exponent partition of the monomial $m$.
Then
$$
\la_i(m^{(k)})-\la_i(m)=\cases
{1, & \hbox{if $1\le i \le k$;} \cr
0, & \hbox{ if $i>k$,} \cr}
$$
and also
$$
\la_i(m')-\la_i(m)\in \{0,1\} \qquad (\forall i);
$$
$$
\sum_{i=1}^n [\la_i(m')-\la_i(m)]=k.
$$
Thus :
$$
\sum_{i=1}^t \la_i(m') \le \sum_{i=1}^t [\la_i(m)+1] =
\sum_{i=1}^t \la_i(m^{(k)}) \qquad (1\le t\le k),
$$
$$
\sum_{i=1}^t \la_i(m') \le k + \sum_{i=1}^t \la_i(m) =
\sum_{i=1}^t \la_i(m^{(k)}) \qquad (k< t\le n).
$$
Therefore
$$
\la(m')\ \underline{\triangleleft}\ \la(m^{(k)}).
$$
If $\la(m')\ {\triangleleft}\ \la(m^{(k)})$ then $m'\prec m^{(k)}$, as claimed.
Otherwise, $\la(m')=\la(m^{(k)})$.

\noindent Now let
$$
I_> \eqdef
\{1\le i \le n\ | \ \lambda_i(m)>\lambda_k(m)\} ;
$$
$$
I_= \eqdef
\{1\le i \le n\ | \ \lambda_i(m)=\lambda_k(m)\} ;
$$
$$
I_< \eqdef
\{1\le i \le n\ | \ \lambda_i(m)<\lambda_k(m)\} .
$$
Namely, $i\in I_>$ if the exponent (in $m$) of $x_{\pi(i)}$ is strictly
larger than that of $x_{\pi(k)}$; and similarly for $I_=$ and $I_<$.
By definition,
$$
\pi(m^{(k)})=\pi(m).
$$
Since $\la(m')=\la(m^{(k)})$, it follows that $\pi(m')$ and $\pi(m)$ agree on 
$I_< \cup I_>$, but may differ on $I_=$. Also, $\pi(m^{(k)})=\pi(m)$ is monotone
increasing on $I_=$. It follows that if $m'\not=m^{(k)}$ then 
$\inv(\pi(m')) > \inv(\pi(m^{(k)}))$ and thus
$m'\prec m^{(k)}$ in the monomial partial order.

\qed

\begin{lem}\label{t.max}
Let $m_1$ and $m_2$ be monomials in $P_n$ of the same total degree,
and let $e_k$ be the $k$-th elementary symmetric function, $1\le k\le n$. 
If $m_1\prec m_2$ then
$$
m_1^{(k)}\prec m_2^{(k)},
$$
where $m_i^{(k)}$ ($i=1,2)$ is the $\prec$-maximal summand in $m_ie_k$  whose existence is 
guaranteed by Lemma 4.2.
\end{lem}

\noindent{\bf Proof.}
By the definition of $m_i^{(k)}$: $\pi(m_i^{(k)})=\pi(m_i)$ and
$\la(m_i^{(k)})=\la(m_i)+\delta^{(k)}$,
where $\delta^{(k)}=(1,\dots,1,0,\dots,0)$ ($k$ ones, $n-k$ zeroes).
Using this with the definition of $\prec$ shows that $m_1\prec m_2
\Longrightarrow m_1^{(k)}\prec m_2^{(k)}$.

\qed

\begin{cor}\label{t.S}
Let $m\in P_n$ be a monomial, $\pi=\pi(m)$ its index permutation, and $\mu=\mu(m)$ 
the complementary partition defined in Subsection 3.2.
Let $S$ be the set of monomials which appear (with nonzero coefficient)
in the expansion of $a_\pi\cdot e_\mu$. Then :
\begin{itemize}
\item[(1)] $m\in S$.
\item[(2)] $m'\in S \hbox{ and } m'\not=m \Longrightarrow m'\prec m$ .
\end{itemize}
\end{cor}

\noindent{\bf Proof.}
Iterative applications of Lemma~\ref{t.mk}, using $\pi(m)=\pi(a_{\pi(m)})$
and Lemma~\ref{t.max}.

\qed

A straightening algorithm follows.

\smallskip

\begin{note}
\noindent{\bf Straightening Algorithm :}

\noindent
For a monomial $m\in P_n$, let $\pi=\pi(m)$ be its index permutation,
$a_\pi$ the corresponding descent basis element, and $\mu=\mu(m)$
the corresponding complementary partition.
Write (by Corollary \ref{t.S})
$$
m=a_\pi\cdot e_\mu - \Sigma,
$$
where $\Sigma$ is a sum of monomials $m'\prec m$. Repeat the process for each $m'$.
\end{note}

\noindent{\bf Example.} Let $m=x_1^2x_2x_3\in P_3$. Then $\pi(m)=123$, $a_{\pi(m)}=1$
 and $\mu(m)=(3,1)$. Hence
$$
x_1^2x_2x_3=1\cdot e_{(3,1)}-\Sigma=(x_1x_2x_3)(x_1+x_2+x_3)-\Sigma,
$$
so $\Sigma=x_1x_2^2x_3+x_1x_2x_3^2$. 
For the first summand $\pi(x_1x_2^2x_3)=213$ and $\mu(x_1x_2^2x_3)=(3)$.
Indeed $x_1x_2^2x_3=a_{213}\cdot e_3=x_2 \cdot (x_1x_2x_3)$.
Similarly, for the second summand 
$\pi(x_1x_2x_3^2)=312$ and $\mu(x_1x_2x_3^2)=(3)$, and
$x_1x_2x_3^2=
a_{312}\cdot e_3=x_3\cdot (x_1x_2x_3)$.
We have obtained
$$
x_1^2x_2x_3=a_{123}\cdot e_{(3,1)}-a_{213}\cdot e_3-a_{312}\cdot e_3.
$$

This algorithm implies

\begin{lem}\label{straight.l}
{\bf (Straightening Lemma)} 
Each monomial $m\in P_n$ has an expression
$$
m= e_{\mu(m)} \cdot a_{\pi(m)}+\sum_{m'\prec m} n_{m',m} e_{\mu(m')} \cdot a_{\pi(m')}, 
$$
where $n_{m',m}$ are integers.
\end{lem}

The following corollary appears in \cite{Al1}.

\begin{cor}\label{4.5'}
The set 
$\{a_\pi+I_n\ |\ \pi\in S_n\}$
forms a basis for the coinvariant algebra $P_n/I_n$.
\end{cor}

\noindent{\bf Proof.}
By Lemma \ref{straight.l}, $\{a_\pi+I_n\ |\ \pi\in S_n\}$
spans $P_n/I_n$ as a vector space over $\bbq$.
Since $\hbox{dim }(P_n/I_n)=|S_n|$ \cite[\S 3.6]{Hum},
this is a basis.

\qed


\subsection{Descent Representations}

Let $\la=(\la_1,\dots,\la_n)$ be a partition with at most $n$ parts. 
Let $P_\la^{\underline{\triangleleft}}$ be the subspace of the polynomial 
ring $P_n=\bbq[x_1,\dots,x_n]$ spanned by
all monomials whose exponent partition is dominated by $\la$:
$$
P_\la^{\underline{\triangleleft}}:=\sn_{\bbq}\{m\ |\ \la(m)\ \underline{\triangleleft}\ \la\}.
$$
Similarly, define $P_\la^{\triangleleft}$ by strict dominance:
$$
P_\la^{\triangleleft}:=\sn_{\bbq}\{m\ |\ \la(m)\triangleleft \la\}.
$$
These subspaces are $S_n$-modules consisting of homogeneous
polynomials of degree $k=|\la|$.

Consider now the canonical projection of $P_n$ onto the coinvariant algebra
$$
\psi : P_n \longrightarrow P_n/I_n.
$$
Define $R_\la$ to be a quotient of images under this map:
$$
R_\la \eqdef \psi(P_\la^{\underline{\triangleleft}})/\psi(P_\la^{\triangleleft}).
$$
Then $R_\la$ is also an $S_n$-module. We will show that every
homogeneous component of the coinvariant algebra may be decomposed into a 
direct sum of certain $R_\la$-s.


For any subset $S\subseteq \{1,\dots,n\}$ define a partition
$$
\lambda_S \eqdef
(\lambda_1,\dots,\lambda_n)
$$
by
$$
\lambda_i \eqdef
|\,S\cap \{i,\dots,n\}\,| \qquad (1\le i\le n).
$$
For example, for any $\pi\in S_n$
$$
\lambda_{\D(\pi)}=(d_1(\pi),\dots,d_n(\pi)).
$$
Hence

\begin{cla}\label{t.g.des}
For any permutation $\pi\in S_n$
$$
\lambda(a_\pi)=\lambda_{\D(\pi)}.
$$
\end{cla}


The Straightening Lemma (Lemma \ref{straight.l}) implies
\begin{lem}\label{t.p}
For any two permutations $\tau$
and $\pi$ in $S_n$, the action of $\tau$ on the monomial $a_\pi
\in  P_n$ has the expression
$$
\tau (a_\pi)=
\sum_{\{w\in S_n\ |\ \lambda(a_w)\ \underline{\triangleleft}\ \lambda(a_\pi)\}}
 n_w a_w + p,
$$
where $n_w\in \bbz$ and $p\in I_n$.
\end{lem}

\noindent{\bf Proof.}
Apply the Straightening Lemma (Lemma \ref{straight.l}) to $m=\tau(a_\pi)$.
Note that $e_{\mu(m')}\not\in I_n$ iff $\mu(m')=\emptyset$, and then
$m'=a_{\pi(m')}$. Denoting $w:=\pi(m')$ we get : 
$m'\preceq m \Longrightarrow \la(a_w)=\la(m')\ \underline{\triangleleft}\ 
\la(m)=\la(\tau(a_\pi))=\la(a_\pi)$.

\qed

Another description of 
the $S_n$-module
$\psi(P^{\underline{\triangleleft}}_\la)$ follows.

\begin{lem}\label{P.la}
For any partition $\la$
$$
\psi(P^{\underline{\triangleleft}}_\la)=
\sn_{\bbq}\{a_\pi+I_n\,|\, \pi\in S_n,\ \la(a_\pi)\ \underline{\triangleleft}\ \la\,\}.
$$
\end{lem}

\noindent{\bf Proof.}
Clearly, for any $\pi\in S_n$ with $\la(a_\pi)\ \underline{\triangleleft}\ \la$,
$a_\pi+I_n=\psi(a_\pi)\in \psi(P^{\underline{\triangleleft}}_\la)$.
The opposite inclusion follows from Lemma~\ref{straight.l}.

\qed

Let
$$
J_\la^{\underline{\triangleleft}} \eqdef
\sn_{\bbq}\{a_\pi+I_n\,|\, \pi\in S_n,\ \la(a_\pi)\ \underline{\triangleleft}\ \la\,\}
$$
and
$$
J_\la^{\triangleleft} \eqdef
\sn_{\bbq}\{a_\pi+I_n\,|\, \pi\in S_n,\ \la(a_\pi)\ {\triangleleft}\ \la\,\}.
$$
By Lemma~\ref{P.la} 
$$
R_{\la}=J_\la^{\underline{\triangleleft}}/J_\la^{\triangleleft}. 
$$
Hence

\begin{cor}\label{R.la}
The following conditions on a partition $\la=(\la_1,\dots,\la_n)$ are equivalent :
\begin{itemize}
\item[(1)] $R_{\la}\ne 0$. 
\item[(2)] $\la=\la(a_\pi)$ for some $\pi\in S_n$.
\item[(3)] $\la=\la_S$ for some $S\subseteq [n-1]$.
\item[(4)] The difference between consecutive parts of $\la$ is either $0$ or $1$,
i.e., (denoting $\la_{n+1}:=0$): 
$$
\la_i-\la_{i+1}\in\{0,1\} \qquad (1\le i\le n).
$$
\end{itemize}
\end{cor}

\noindent{\bf Proof.}
$R_{\la}\ne 0$ if and only if $J_\la^{\underline{\triangleleft}}\ne J_\la^{\triangleleft}$. 
This happens if and only if there exists a permutation $\pi\in S_n$ for which
$\la(a_\pi)=\la$ (since, by Corollary~\ref{4.5'}, the various $a_{\pi}+I_n$ are linearly independent).
By Claim \ref{t.g.des}, $\la(a_\pi)=\la_S$ for $S=\D(\pi)$.
Conversely, any subset of $[n-1]$ is a descent set for an appropriate permutation. This shows the equivalence of the first three conditions.
The equivalence of (3) and (4) follows easily from the definition of $\la_S$.

\qed

From now on denote $R_S:=R_{\la_S}$.
For $\pi\in S_n$, let $\bar a_\pi$ be the image of the descent basis element
$a_\pi+I_n\in J_{\la_S}^{\underline{\triangleleft}}$ in the quotient
$R_{S}$, where $S:=\D(\pi)$.

\begin{cor}\label{4.8}
For any subset $S\subseteq [n-1]$, the set 
$$
\{\,\bar{a}_\pi\,|\, \pi\in S_n,\ \D(\pi)=S\,\}
$$
forms a basis of $R_{S}$.
\end{cor}

\noindent{\bf Proof.}
Follows by elementary linear algebra from Corollary \ref{4.5'},
Claim \ref{t.g.des},
and the definitions of 
$J_S^{\underline{\triangleleft}}$, $J_S^{\triangleleft}$ and $\bar{a}_\pi$.

\qed

Recall the notation $R_k$ for the $k$-th homogeneous component
of the coinvariant algebra $P_n/I_n$.

\begin{thm}\label{a.main}
For every $0\le k\le {n\choose 2}$, 
$$
R_k \cong \bigoplus_{S} R_S
$$
as $S_n$-modules, where the sum is over all subsets
$S\subseteq [n-1]$ such that $\sum_{i\in S} i =k$.
\end{thm}

\noindent{\bf Proof.}
Note that $a_\pi+I_n\in R_k \Longleftrightarrow \maj(\pi)=k$.
It follows, by Corollary \ref{4.5'}, that
$\{a_\pi+I_n|\maj(\pi)=k\}$ is a basis for $R_k$;
so that, by Corollary \ref{4.8},
$R_k \cong \bigoplus_{S} R_S$
(sum  over all $S\subseteq [n-1]$ with $\sum_{i\in S} i =k$)
as vector spaces over $\bbq$.
By Maschke's Theorem, if $V$ is a finite dimensional $G$-module
for a finite group $G$ (over a field of characteristic zero)
and $W\subseteq V$ is a $G$-submodule, then $V\cong W\oplus (V/W)$ as $G$-modules.
Apply this to the poset 
$$
\{J_S^{\underline{\triangleleft}}\ |\ S\subseteq [n-1]\ ,\ \sum_{i\in S} i =k\}
$$
ordered by dominance order on the partitions $\la_S$.
Using the definition of $R_S$, we get by induction on the poset
the required isomorphism.

\qed

\section{Decomposition of Descent Representations}

In this section we refine the Lusztig-Stanley Theorem.

\begin{thm}\label{t.rep.decomp}
For any subset $S\subseteq [n-1]$ and partition $\mu\vdash n$,
the multiplicity in $R_S$ of the irreducible $S_n$-representation
corresponding to $\mu$ is
$$
m_{S,\mu} \eqdef
|\,\{\,T\in SYT(\mu)\ |\ \D(T)=S\,\}\,| .
$$
\end{thm}

The proof applies an argument of Stanley (cf. \cite[Theorem 8.8]{Re}
and references therein), using symmetric functions. 
Its use is possible here due to the Straightening Lemma.
A key lemma in the proof is a multivariate version of Stanley's formula.

\subsection{A Multivariate Version of Stanley's Formula}

In this subsection we state and prove a multivariate version of a well known 
formula of Stanley for the principal specialization of a Schur function.

Recall that $s_{\la}(x_1,\ldots,x_n)$ denotes the Schur function corresponding 
to a partition $\la$ of $n$.  Then

\begin{pro}\label{t.stanley}
\cite[Prop. 7.19.11]{StEC2}
If $\la$ is a partition of $n$ then 
$$
s_{\la}(1,q,q^2,\dots,q^{n-1})=
{\sum_{T\in SYT(\la)} q^{\maj(T)} \over (1-q)(1-q^2)\cdots(1-q^n)},
$$
where $T$ runs through all standard Young tableaux of shape $\la$.
\end{pro}


Define $\bbq[[z_1,z_1z_2,\ldots]]$ to be the ring of formal power series 
in countably many variables $z_1,z_1z_2,\dots,z_1z_2\cdots z_k,\dots$; 
a linear basis for it consists of the monomials 
$z^\la:=z_1^{\la_1}\cdots z_n^{\la_k}$ 
for all partitions $\la=(\la_1,\dots,\la_k)$ ($\la_1\ge \dots\ge \la_k\ge 0$). 
Let $\bbq[[q_1,q_1q_2,\ldots]]$ be similarly defined.

Let $\iota:\bbq[[z_1,z_1z_2,\ldots]]\longrightarrow \bbq[[q_1,q_1q_2,\ldots]]$ 
be defined by
$$
\iota(z^\la)\eqdef
q^{\la'} \qquad (\forall \la), \leqno(4.1)
$$
where $\la'$ is the partition conjugate to $\la$.
Extend $\iota$ by linearity. Note that $\iota$ is {\em not} a ring homomorphism.


For any standard Young tableau $T$ define
$$
d_i(T) \eqdef |\,\{\,j\ge i\,|\,j\in \D(T)\,\}\,| \qquad(1\le i\le n),
$$
where $\D(T)$ is the set of all descents in $T$ 
(as in Subsection~\ref{s.part_tab}).

\smallskip

\begin{lem}\label{t.tau}{\bf (Multivariate Formula)} 
If $\la$ is a partition of $n$ then
$$
\iota[s_\la(1,z_1,z_1 z_2,\ldots,z_1 z_2\cdots z_i,\dots)]=
{\sum_{T\in SYT(\la)} \prod_{i=1}^n q_i^{d_i(T)}
\over \prod_{i=1}^n (1 - q_1 q_2 \cdots q_i)},
$$
where $T$ runs through all standard Young tableaux of shape $\lambda$.
\end{lem}

\smallskip

Note that Proposition~\ref{t.stanley} is the special case obtained by
substituting $q_1=q_2=\ldots=q_n=q$.
Lemma~\ref{t.tau} is actually implicit in~\cite[p.~27]{St72} 
(formula (18), within the proof of Theorem~9.1), in the general context of 
$(P,\omega)$-partitions.
The fact that the LHS there is equal to the LHS above is part of the 
following proof, which will later be adapted to prove a $B$-analogue 
(Lemma~\ref{t.tau.B} below).
The key bijection is basically the one used in~\cite[Lemma~3.1]{Stem}.

\smallskip

\noindent{\bf Proof.}
By Claim~\ref{t.schur} 
$$
s_{\la}(x_1,x_2,\dots)=
\sum_{\hT\in RSSYT(\la)} 
\prod_{i=1}^\infty x_i^{m_i(\hT)},
$$
where $\hT$ runs through all reverse semi-standard Young tableaux of 
shape $\la$; and
$$
m_i(\hT) \eqdef
|\{\hbox{cells in $\hT$ with entry $i$}\}| \qquad (\forall i\ge 1).
$$
Letting
$$
x_1=1
$$
and
$$
x_i = z_1 z_2 \cdots z_{i-1} \qquad (i\ge 2)
$$
we get
$$
s_\lambda(1,z_1,z_1 z_2,\dots)=
\sum_{\hT\in RSSYT(\la)}  
\prod_{i=1}^\infty z_i^{m_{> i}(\hT)},
$$
where
$$
m_{> i}(\hT) \eqdef
|\{\hbox{cells in $\hT$ with entry $> i$}\}| \qquad (\forall i).
$$
Of course, 
$$
m_{> 0}(\hat T)=n.
$$
Let $\mu(\hT)$ be the vector $(m_{> 1}(\hT),m_{> 2}(\hT),\dots)$.
Then $\mu(\hT)$ is a partition with largest part at most $n$.
Note that the conjugate partition is
$$
\mu(\hT)'=(\hT_1 -1,\hT_2 -1,\dots,\hT_n -1),
$$
where $\hT_1, \ldots, \hT_n$ are the entries of $\hT$ 
in weakly decreasing order.
Thus
$$
\iota[s_{\la}(1,z_1,z_1 z_2,\dots)]=
\sum_{\hT\in RSSYT(\la)}  
\prod_{i=1}^n q_i^{\hT_i -1}. \leqno (4.2)
$$

\smallskip

Recall that $SYT(\la)$ and $RSSYT(\la)$ 
are the sets of standard Young tableaux and 
reverse semi-standard Young tableaux of shape $\la$, respectively.
For any $\la\vdash n$ define a map
$$
\phi_{\la}\,:\,RSSYT(\la) \longrightarrow SYT(\la) \times \bbn^n 
$$
by
$$
\phi_{\la}(\hT) \eqdef (T,\Delta),
$$
where, for $\hT \in RSSYT(\la)$, $T$ is a standard Young tableau of shape $\la$ and
$\Delta=(\Delta_1,\dots,\Delta_n)$ is a sequence of nonnegative integers,
defined as follows:
\begin{itemize}
\item[(1)] 
Let $(\hT_1,\ldots,\hT_n)$ be the vector of entries of $\hT$, in weakly decreasing 
order.  Then $T$ is the standard Young tableau, of the same shape as $\hT$, 
having entry
$i$ ($1\le i\le n$) in the same cell in which $\hT$ has entry $\hT_i$. 
If some of the  entries of $\hT$ are equal then they necessarily belong to distinct
columns, and the corresponding entries of $T$ are then chosen 
increasing from left to right (i.e., with increasing column indices).
\item[(2)]
Define
$$
\Delta_i\eqdef \hT_i - d_{i}(T) - \hT_{i+1} + d_{i+1}(T) \qquad(1\le i\le n),
$$
where, by convention, $\hT_{n+1} \eqdef 1$ and $d_{n+1}(T) \eqdef 0$.
%
%
\end{itemize}

\smallskip

\noindent{\bf Example.}
Let $\la=(3,2,2)$ and
$$
\hT=
\begin{array}{ccc}
7 & 4 & 4 \\
4 & 2 &   \\
3 & 1 &   \\
\end{array}
\in RSSYT(\la).
$$
Computing  $\phi_{\la}(\hT) = (T,\Delta)$, the first step yields
$$
T=
\begin{array}{ccc}
1 & 3 & 4 \\
2 & 6 &   \\
5 & 7 &   \\
\end{array}
\in SYT(\la),
$$
so that $\D(T)=\{1,4,6\}$.  
Therefore, 
$(\hT_1-d_1(T),\dots,\hT_7-d_7(T))=(4,2,2,2,2,1,1)$ and
$\Delta=(2,0,0,0,1,0,0)$.

\smallskip

The following claim is easy to verify.
\begin{cla}\label{t.phi}

\begin{itemize}
\item[{(1)}] $\phi_\la$ is a bijection.
\item[{(2)}] If $(\hT_1,\dots,\hT_n)$ is the vector of entries of $\hT$,
in weakly decreasing order, and $\phi_{\la}(\hT)=(T,\Delta)$, then
$$
\hT_i -1 = d_i(T)+\sum_{j \geq i} \Delta_j \qquad(1\le i\le n).
$$
\end{itemize}
\end{cla}
By Claim~\ref{t.phi}, for any $\hT\in RSSYT(\la)$
$$
\prod_{i=1}^n q_i^{\hat T_i -1}=\prod_{i=1}^n q_i^{d_i(T)}\cdot
\prod_{j=1}^n (q_1\cdots q_j)^{\Delta_j}, \leqno(4.3)
$$
where $(T,\Delta)=\phi_{\la}(\hT)$.

Substituting (4.3) into (4.2) we get the formula in the statement of
Lemma~\ref{t.tau}.

\qed

\subsection{Proof of Theorem~\ref{t.rep.decomp}}

For a permutation $\tau\in S_n$
let the (graded) {\it trace} of its action on the polynomial ring $P_n$ be
$$
\tr_{P_n}(\tau) \eqdef
\sum_{m} \langle \tau(m),m\rangle \cdot {\bar q}^{\la(m)},
$$
where the sum is over all monomials $m$ in $P_n$, $\lambda(m)$
is the exponent partition of the monomial $m$, and the 
inner product is such that the set of all monomials is an orthonormal
basis for $P_n$.  In particular, $\langle \tau(m),m\rangle\in \{0,1\}$
$(\forall \tau\in S_n)$.

\begin{cla}\label{t.tr}
If $\tau\in S_n$ is of cycle type $\mu$ then the trace of its action on $P_n$ is
$$
\tr_{P_n}(\tau)=\iota[p_\mu(1,z_1,z_1 z_2,\dots)],
$$
where $p_\mu$ is the power sum symmetric function corresponding to $\mu$, and $z_1,z_2,\ldots$ are independent variables.
\end{cla}

\noindent{\bf Proof.}
For a monomial $m$ in $P_n$, $\langle\tau(m),m\rangle = 1$ (i.e., $\tau(m) = m$)
if and only if, for each cycle of $\tau$, the variables $x_i$ with indices in that cycle 
all have the same exponent in $m$.  If $\mu = (\mu_1,\ldots,\mu_t)$ is the cycle type
of $\tau$ we thus obtain, for each sequence $(e_1,\ldots,e_t)\in \bbn^t$, a unique 
monomial $m$ with $\tau(m) = m$ and $e_j$ as the common exponent for cycle 
number $j$ ($1\le j\le t$).   Thus $\la(m)$ consists of $\mu_j$ copies of $e_j$ ($1\le j\le t$),
reordered; and 
$$
\iota^{-1}[{\bar q}^{\la(m)}] = {\bar z}^{\la(m)'} = 
\prod_{j=1}^{t} (z_1 \cdots z_{e_j})^{\mu_j}.
$$
Summing over all choices of $e_1,\ldots,e_t$ gives $p_{\mu}(1,z_1,z_1 z_2,\ldots)$.

\qed

Recall that for the coinvariant algebra $P_n/I_n$ we have
the descent basis $\{a_\pi+I_n\ |\ \pi\in S_n\}$.
Define, for $\tau\in S_n$:
$$
\tr_{P_n/I_n}(\tau) \eqdef
\sum_{\pi\in S_n} \langle \tau(a_\pi+I_n),a_\pi+I_n\rangle \cdot  
{\bar q}^{\lambda(a_\pi)},
$$
where the inner product is such that the descent basis is orthonormal.

From the Straightening Lemma (Lemma \ref{straight.l})
it follows that
\begin{cla}\label{t.trtr}
For every $n\ge 1$ and every $\tau\in S_n$
$$
\tr_{P_n}(\tau)=\tr_{P_n/I_n}(\tau)\cdot \sum_{\lambda} \bar q^\lambda,
$$
where the sum is over all partitions $\lambda$ with at most $n$ parts.
\end{cla}

\noindent{\bf Proof.}
Replace the monomial basis of $P_n$ by the homogeneous polynomial basis
$\{a_\pi e_\mu|\ \pi\in S_n,\ \mu \hbox{ is a partition with largest part at most } n\}$.
The trace $\tr_{P_n}(\tau)$ is not changed, provided that we now use the inner product
for which the new basis is orthonormal. 
Note that $\tau(e_\mu)=e_\mu$ ($\forall\tau\in S_n$).
Also note that, by the Straightening Lemma, the $\prec$-maximal
monomial $m$ in $a_\pi e_\mu$ has an exponent partition
$\la(m)=\la(a_\pi)+\mu'$, where addition of partitions is componentwise. 

\qed

\begin{obs}\label{t.lambda}
$$
\sum_{\la} \bar q^\la = \prod_{i=1}^n{1\over 1-q_1\cdots q_i}\,,
$$
where the sum is over all partitions $\lambda$ with at most $n$ parts.
\end{obs}

Using the Frobenius formula, the multivariate version of Stanley's formula (Lemma~\ref{t.tau}) and Claim~\ref{t.tr}, 
we obtain for every permutation $\tau\in S_n$ of cycle type $\mu$:
\begin{eqnarray*}
\tr_{P_n}(\tau) 
&=& \iota[p_\mu(1,z_1,z_1 z_2,\dots)]=
    \sum_{\la\vdash n}\chi^\la_\mu \iota[s_\la(1,z_1,z_1z_2,\ldots)]=\\
&=& \sum_{\la\vdash n}\chi^\lambda_\mu \cdot
       {\sum_{T\in SYT(\la)} \prod_{i=1}^n q_i^{d_i(T)}
        \over \prod_{i=1}^n (1 - q_1 q_2 \cdots q_i)}.
\end{eqnarray*}
By Claim~\ref{t.trtr} and Observation~\ref{t.lambda} we now get
$$
\tr_{P_n/I_n}(\tau)=
\sum_{\la\vdash n}\chi^\la_\mu  
\sum_{T\in SYT(\la)}
\prod_{i=1}^n q_i^{d_i(T)}.
$$
We conclude that the graded multiplicity in $P_n/I_n$ of the irreducible $S_n$-representation 
corresponding to $\la$ is
$$
\sum_{T\in SYT(\la)} \prod_{i=1}^{n} q_i^{d_i(T)} = \sum_{T\in SYT(\la)} {\bar q}^{\la_{\D(T)}}.
$$
Consider now the claim of Theorem~\ref{a.main}.  Its proof shows that 
$$
P_n/I_n \cong \bigoplus_{S \subseteq [n-1]} R_S
$$
actually holds as an isomorphism of {\em graded} $S_n$-modules.  By Corollary~\ref{4.8}
and Claim~\ref{t.g.des}, $R_S$ is the homogeneous component of multi-degree $\la_S$ in
$P_n/I_n$.  It thus follows that the multiplicity in $R_S$ of the irreducible $S_n$-representation
corresponding to $\la$ is 
$$
|\,\{T\in SYT(\la)\,|\,\D(T) = S\}\,|,
$$
and the proof of Theorem~\ref{t.rep.decomp} is complete.

\qed

\section{Descent Representations for Type $B$}

In this section we give $B$-analogues of the concepts and results in Sections 3 and 4.

\subsection{The Signed Descent Basis}

For any signed permutation $\sigma\in B_n$ let
\[ 
\D(\sigma ) \eqdef \{ i \in [n-1]: \; \sigma (i) > \sigma(i+1) \}  
\]
be the set of descents in $\sigma$ with respect to the standard linear order on the integers,
and let
\[
d_{i}(\sigma ) \eqdef |\{ j \in \D(\sigma ): \; j \geq i \} | \qquad (1\le i\le n)
\]
be the number of descents in $\sigma$ from position $i$ on.

Also let
$$
 \varepsilon _{i} (\sigma ) \eqdef \left\{
\begin{array}{ll}  
1, & \mbox{if $\sigma  (i)<0$,} \\
0, & \mbox{otherwise,}
\end{array} \right. 
$$
and
$$
f_i(\sigma):=2d_i(\sigma)+\varepsilon_i(\sigma).
$$

Note that 
$$
\sum_{i=1}^n f_i(\sigma)=\fm(\sigma).
$$

\smallskip

Finally, associate to $\sigma$ the monomial
$$
b_\sigma \eqdef
\prod_{i=1}^n x_{|\sigma(i)|}^{f_i(\sigma)}.
$$


\smallskip

We will show that the set $\{{b_\sigma}+I^B_n\,|\,\sigma\in B_n\}$ forms
a linear basis for the coinvariant algebra of type $B$.
We call it the {\it signed descent basis}.

\subsection{Straightening}

The {\it signed index permutation} of a monomial $m=\prod_{i=1}^n x_i^{p_i}\in P_n$ is 
the unique signed permutation $\sigma=\sigma(m) \in B_n$ such that
$$
p_{|\sigma(i)|}\ge p_{|\sigma(i+1)|} \qquad (1\le i < n) , \leqno(1)
$$
$$
p_{|\sigma(i)|}=p_{|\sigma(i+1)|}  \Longrightarrow
 \sigma(i)<\sigma(i+1) \leqno(2)
$$
and
$$
p_{|\sigma(i)|} \equiv 0 (\hbox{mod } 2) \Longleftrightarrow  \sigma(i)>0 . \leqno(3)
$$
In other words, $\sigma$ reorders the variables $x_i$ as does the corresponding 
index permutation of type A (see Subsection~\ref{s.index_exp}), after attaching
a minus sign to (indices of) variables with odd exponents.  
Note that this reverses the order of ``negative'' indices.

\smallskip

\noindent{\bf Example.} 
 $m=x_1^2 x_2^3 x_3^2 x_5 x_6^3$ and $n=7$.
Then $m=x_2^3 x_6^3 x_1^2 x_3^2 x_5^1 x_4^0 x_7^0$ and 
$\sigma(m)=[-6,-2,1,3,-5,4,7]\in B_7$.

\smallskip

Let $m=\prod_{i=1}^n x_i^{p_i}$ be a monomial in $P_n$,
$\sigma=\sigma(m)$ its signed index permutation,
and $b_\sigma$ the corresponding signed descent basis element. 
In analogy with Claim \ref{t.exp}
(for type $A$) we have
\begin{cla}\label{b.1}
The sequence $(p_{|\sigma(i)|}-f_i(\sigma))_{i=1}^n$
of exponents in ${m/ b_\sigma}$ consists of nonnegative even integers, and
is weakly decreasing:
$$
p_{|\sigma(i)|}-f_i(\sigma)\ge 
p_{|\sigma(i+1)|}-f_{i+1}(\sigma) \qquad (1\le i< n) .
$$
\end{cla}

\noindent{\bf Proof.}
By condition (3) above 
$\varepsilon_i(\sigma)=0\Longleftrightarrow 
p_{|\sigma(i)|}\equiv 0 (\hbox{mod }2)$.
The numbers
$p_{|\sigma(i)|}-f_i=p_{|\sigma(i)|}-\varepsilon_i(\sigma)-2d_i(\sigma)$
are therefore even integers. Also $d_n(\sigma)=0$ so that
$p_{|\sigma(n)|}-f_n=p_{|\sigma(n)|}-\varepsilon_n(\sigma)\ge 0$.
It remains to show that the sequence is weakly decreasing.

If $\sigma(i)<\sigma(i+1)$ then $d_i(\sigma)=d_{i+1}(\sigma)$, 
and $p_{|\sigma(i+1)|}-\varepsilon_{i+1}(\sigma)\ge
p_{|\sigma(i)|}-\varepsilon_i(\sigma)$ since their difference is an even integer $\ge -1$.

If $\sigma(i)>\sigma(i+1)$ 
and $\varepsilon_i(\sigma)=\varepsilon_{i+1}(\sigma)$
then $d_i(\sigma)=d_{i+1}(\sigma)+1$;
so $(p_{|\sigma(i+1)|}-f_{i+1}(\sigma))-
(p_{|\sigma(i)|}-f_{i}(\sigma))=p_{|\sigma(i+1)|}-p_{|\sigma(i)|}-2\ge 0$ \ since\ $p_{|\sigma(i+1)|}\equiv p_{|\sigma(i+1)|}(\hbox{mod }2)$\ 
and because of conditions (1) and (2) above.

Finally, if $\sigma(i)>\sigma(i+1)$ 
and $\varepsilon_i(\sigma)\not=\varepsilon_{i+1}(\sigma)$
then $\sigma(i)>0>\sigma(i+1)$ so that 
$\varepsilon_i(\sigma)=0$, $\varepsilon(i+1)=1$ and
$d_i(\sigma)=d_{i+1}(\sigma)+1$. Thus
$(p_{|\sigma(i+1)|}-f_{i+1}(\sigma))-
(p_{|\sigma(i)|}-f_{i}(\sigma))=p_{|\sigma(i+1)|}-p_{|\sigma(i)|}-1\ge 0$.

\qed

Denote by $\mu_B(m)$
the partition conjugate to 
$({p_{|\sigma(i)|}-f_i(\sigma)\over 2})_{i=1}^n$ .

\smallskip

\noindent{\bf Definition:} For monomials $m_1,m_2\in P_n$ of the same total degree,
$m_1\prec_B m_2$ if, for every $1\le i\le n$, the exponents of $x_i$ in $m_1$ and $m_2$
have the same parity; and also either
\begin{itemize}
\item[(1)] $\lambda(m_1)\ \triangleleft\ \lambda(m_2)$ (strictly smaller in dominance order); or
\item[(2)] $\lambda(m_1) = \lambda(m_2)$ and $\inv (\pi(m_1)) > \inv (\pi(m_2))$.
\end{itemize}
Here $\pi(m_i)$ is the (unsigned) index permutation
of $m_i$ as in Subsection 3.2.

\smallskip

Imitating the proof of Lemma~\ref{straight.l} we obtain

\begin{cor}\label{b.2} {\bf (Straightening Lemma)}
Each monomial $m\in P_n$
has an expression
$$
m= e_{\mu_B(m)}(x_1^2,\dots,x_n^2)\cdot b_{\sigma(m)} + \sum_{m'\prec_B m} 
n_{m',m}e_{\mu_B(m')}(x_1^2,\dots,x_n^2) \cdot b_{\sigma(m')},
$$
where $n_{m',m}$ are integers.
\end{cor}

\subsection{Construction of Descent Representations}

Recall (from Subsection 2.5) the definition of $I^B_n$.

\begin{cor}\label{6.2'}
The set 
$$
\{b_\sigma+I^B_n\ |\ \sigma\in B_n\}
$$
forms a basis for the coinvariant algebra $P_n/I^B_n$.
\end{cor}

\noindent {\bf Proof.} Similar to the proof of Corollary \ref{4.5'}.

\qed

\noindent{\bf Remark.} A different proof of Corollary \ref{6.2'}
may be obtained by combining and extending the proofs of 
\cite[Theorems 1(3) and 2]{ATY}.

\begin{cor}\label{b.3}
For any pair of elements $\tau$
and $\sigma$ in $B_n$, the action of $\tau$ on $b_\sigma\in P_n$ 
has the expression 
$$
\tau (b_\sigma)=
\sum_{\{w\in B_n\ |\ \la(b_w)\ \underline{\triangleleft}\ \la(b_\sigma)\}} 
n_w b_w + p,
$$
where $n_w\in \bbz$ and $p\in I_n^B$.
\end{cor}

\noindent{\bf Proof.} The same as the proof of Lemma \ref{t.p}.
\qed

\smallskip

Note that, actually, $n_w\not=0$ in Corollary \ref{b.3}
only if $b_w$ and $b_\sigma$ have the same number of odd exponents.
%
Note also that, unlike
the case of type $A$ (Claim~\ref{t.g.des}), for $\sigma\in B_n$,
$\la(b_\sigma)$ and $\la_{\D(\sigma)}$ are not necessarily equal.
Indeed, for subsets $S_1\subseteq [n-1]$
and $S_2\subseteq [n]$ define a vector $\la_{S_1,S_2}$ by
$$
\la_{S_1,S_2}:=2\la_{S_1}+{\mathbf 1}_{S_2},
$$
where $\la_{S_1}$ is as in Subsection 3.5 above, ${\mathbf 1}_{S_2}\in
\{0,1\}^n$ is the characteristic vector of $S_2$, and addition is componentwise. Thus
$$
\la_{S_1,S_2}(i)=2\cdot |\{j\ge i|\ j\in S_1\}|+{\mathbf 1}_{S_2}(i).
$$
It should be noted that $\la_{S_1,S_2}$ is not always a partition.

\begin{cla}\label{5.4'}
For any $\sigma\in B_n$
$$
\la(b_\sigma)
=\la_{S_1,S_2}
$$ 
where $\D(\sigma)=S_1$ and $\N(\sigma)=S_2$.
\end{cla}

Now define
$$
R^B_{\la}:=\psi^B(P^{\underline{\triangleleft}}_{\la})\ /\ 
\psi^B(P^{\triangleleft}_{\la}),
$$
where $\psi^B:P_n\longrightarrow P_n/I^B_n$ is the canonical map
from $P_n$ onto the coinvariant algebra of type $B$,
and $P^{\underline{\triangleleft}}$, $P^{\triangleleft}$
are as in Subsection 3.5.

By arguments similar to those given in the proof of Corollary
\ref{P.la}, $R^B_\la$ may be described via
the signed descent basis :
$$
R^B_\la=J^{\underline{\triangleleft}}_{\la,B} /J^{\triangleleft}_{\la,B}, 
$$
where
$$
J^{\underline{\triangleleft}}_{\la,B} \eqdef
\sn_{\bbq}\{b_\sigma+I^B_n\ |\ 
\sigma\in B_n,\ 
\la(b_\sigma)\ {\underline{\triangleleft}}\ \la\},
$$
and
$$
J^{\triangleleft}_{\la,B} \eqdef
\sn_{\bbq}\{b_\sigma+I^B_n\ |\ 
\sigma\in B_n,\ 
\la(b_\sigma)\ {\triangleleft}\ \la\}.
$$

\begin{cor}\label{R.la.b}
The following conditions on a partition $\la=(\la_1,\dots,\la_n)$ are equivalent :
\begin{itemize}
\item[(1)] $R^B_{\la}\ne 0$. 
\item[(2)] $\la=\la(b_\sigma)$ for some $\sigma\in B_n$.
\item[(3)] $\la=\la_{S_1,S_2}$ for some $S_1\subseteq [n-1]$ and 
$S_2\subseteq [n]$, and $\la_{S_1,S_2}$ is a partition.
\item[(4)] The difference between consecutive parts of $\la$ is either $0$, $1$ or $2$,
i.e. (denoting $\la_{n+1}:=0$): 
$$
\la_i-\la_{i+1}\in\{0,1,2\} \qquad (1\le i\le n).
$$
\end{itemize}
\end{cor}

\noindent{\bf Proof.} Similar to the proof of Corollary~\ref{R.la}.
\qed

\smallskip

From now on, denote $R^B_{S_1,S_2}:=R^B_{\la_{S_1,S_2}}$.
For $\sigma\in B_n$,
denote by $\bar b_\sigma$ the image of the signed descent basis element
$b_\sigma+I^B_w\in J^{\underline{\triangleleft}}_{\la_{S_1,S_2}, B}$ in 
$R^B_{S_1,S_2}$.

\begin{cor}\label{4.8.b}
For any $S_1\subseteq [n-1]$ and $S_2\subseteq [n]$, the set  
$$
\{\,\bar{b}_\sigma\,|\, \sigma\in B_n,\ \D(\sigma)=S_1,\ \N(\sigma)=S_2\}
$$
forms a basis of $R^B_{S_1,S_2}$.
\end{cor}

Recall the notation $R^B_k$ for the $k$-th homogeneous component
of $P_n/I^B_n$.

\begin{thm}\label{b.4}
For every $0\le k\le n^2$,
$$
R^B_k \cong \bigoplus_{S_1,S_2} R^B_{S_1,S_2},
$$
as $B_n$-modules, where the sum is over all subsets
$S_1\subseteq [n-1]$ and $S_2\subseteq [n]$ such that
$\lambda_{S_1,S_2}$ is a partition and
$$
2\cdot \sum_{i\in S_1} i +|S_2|=k.
$$
\end{thm}

\noindent{\bf Proof.}
Similar to the proof of Theorem \ref{a.main}.
\qed

\subsection{Decomposition of Descent Representations}

In this subsection we give a $B$-analogue of Theorem~\ref{t.rep.decomp}.

\begin{thm}\label{t.rep.decomp.B}
For any pair of subsets $S_1\subseteq [n-1]$, $S_2\subseteq [n]$,
and a bipartition $(\mu^1,\mu^2)$ of $n$, the multiplicity of
the irreducible $B_n$-representation corresponding to $(\mu^1,\mu^2)$
in $R^B_{S_1,S_2}$ is
$$
m_{S_1,S_2,\mu^1,\mu^2} \eqdef
|\,\{\,T\in SYT(\mu^1,\mu^2)\ | \; \D(T)=S_1 \hbox{ and } \N(T)=S_2\,\}\,| .
$$
\end{thm}

Again, we need a (type $B$) multivariate version of Proposition~\ref{t.stanley}.
Recall definition~(4.1) of the mapping $\iota$.

\smallskip

Let $\lala$ be a bipartition of $n$.
Recall the definitions of a standard Young tableau $T=\TT$ of shape $\lala$  
and the sets $\D(T)$ and $\N(T)$ 
from Subsection 2.3. Denote the set of all standard tableaux of shape
$\lala$ by  $SYT\lala$.
A {\em reverse semi-standard Young tableau of shape $\lala$} is a pair
$\hTT$ of reverse semi-standard Young tableaux, where:
$\hT^i$ has shape $\la^i$ $(i=1,2)$; the entries of $\hT^i$ are congruent to $i$ 
$(\mbox{\rm mod\ } 2)$; the entries in each row are weakly decreasing; and the entries
in each column are strictly decreasing.
Denote by $RSSYT\lala$ the set of all such tableaux.

\smallskip

For any standard Young tableau $T=\TT$ of shape $\lala$ define
$$
d_i(T) \eqdef |\{j\in \D(T) \,|\, j\ge i\}|,
$$
$$
 \varepsilon _{i} (T) \eqdef \left\{
\begin{array}{ll}  
1, & \mbox{if $i\in \N(T)$,} \\
0, & \mbox{otherwise,}
\end{array} \right. 
$$
and
$$
f_i(T):=2 d_i(T)+\varepsilon_i(T).
$$
Then
$$
\fm(T)=\sum_{i=1}^n f_i(T).
$$

\begin{lem}\label{t.tau.B}
If $\lala$ is a bipartition of $n$ then
$$
\iota[s_{\la^1}(1, z_1 z_2, z_1 z_2 z_3 z_4,\ldots) \cdot
     s_{\la^2}(z_1, z_1 z_2 z_3,\ldots)] =
$$
$$=
{\sum_{T\in SYT\lala} \prod_{i=1}^n q_i^{f_i(T)}
\over \prod_{i=1}^n (1 - q_1^2 q_2^2 \cdots q_i^2)},
$$
where $T$ runs through all standard Young tableaux of shape $\lala$.
\end{lem}

\smallskip

\noindent{\bf Proof.}
By Claim~\ref{t.schur} 
$$
s_{\la}(x_1,x_2,\dots)=
\sum_{\hT\in RSSYT(\la)} \prod_{i=1}^\infty x_i^{m_i(\hT)},
$$
where $\hT$ runs through all reverse semi-standard Young tableaux of shape $\la$, and
$$
m_i(\hT) \eqdef
|\{\hbox{cells in $\hT$ with entry $i$}\}| \qquad (\forall i\ge 1).
$$
Note that for $(\hT^1,\hT^2)\in RSSYT(\la^1,\la^2)$ we have 
$(\hT^1 +1)/2\in RSSYT(\la^1)$, where  $(\hT^1 +1)/2$ is obtained 
by replacing each entry $i$ of $\hT^1$ by $(i+1)/2$;
and similarly $\hT^2/2\in RSSYT(\la^2)$.
Thus
$$
s_{\la^1}(x_1,x_2,\dots)\cdot s_{\la^2}(y_1,y_2,\dots)=
\sum_{\hTT\in RSSYT\lala}  
\prod_{i=1}^\infty x_i^{m_{2i-1}(\hT^1)}y_i^{m_{2i}(\hT^2)}.
$$
Letting
$$
x_1=1 
$$
and
$$
x_i = z_1 z_2 \cdots z_{2i-2} \qquad (i\ge 2)
$$
for $s_{\la^1}$, while
$$
y_i = z_1 z_2 \cdots z_{2i-1} \qquad (i\ge 1)
$$
for $s_{\la^2}$, we get
$$
s_{\la^1}(1, z_1 z_2, \ldots) \cdot s_{\la^2}(z_1, z_1 z_2 z_3, \ldots) =
\sum_{\hTT\in RSSYT\lala}  
\prod_{i=1}^\infty z_i^{m_{> i}\hTT},
$$
where
$$
m_{> i}\hTT \eqdef
|\{\hbox{cells in $\hTT$ with entry $> i$}\}| \qquad (\forall i\ge 1).
$$
Let $\mu\hTT$ be the vector $(m_{> 1}\hTT,m_{> 2}\hTT,\ldots)$.
Then $\mu\hTT$ is a partition with largest part at most $n$.
The conjugate partition is
$$
\mu\hTT'=(\hTT_1 -1,\hTT_2 -1,\dots,\hTT_n -1),
$$
where $\hTT_1,\dots,\hTT_n$ are the entries of $\hTT$ 
in weakly decreasing order.

Thus
$$
\iota[s_{\la^1}(1,z_1 z_2,\dots) \cdot s_{\la^2}(z_1, \ldots)]=
\sum_{\hTT\in RSSYT\lala}  
\prod_{i=1}^n q_i^{\hTT_i -1}. \leqno (5.1)
$$

\smallskip

For any bipartition $\lala$ of $n$ define a map
$$
\phi_{\lala}\,:\,RSSYT\lala \longrightarrow SYT\lala \times \bbn^n 
$$
by
$$
\phi_{\lala}\hTT \eqdef (\TT,\Delta),
$$
where, for $\hTT \in RSSYT\lala$, $\TT$ is a standard Young tableau 
of shape $\lala$ and
$\Delta=(\Delta_1,\dots,\Delta_n)$ is a sequence of nonnegative integers,
defined as follows:
\begin{itemize}
\item[(1)] 
Let $(\hTT_1,\ldots,\hTT_n)$ be the vector of entries of $\hTT$, in weakly decreasing 
order.  Then $\TT$ is the standard Young tableau, of the same shape as $\hTT$, with
entry $i$ ($1\le i\le n$) in the same cell in which $\hTT$ has entry $\hTT_i$. 
If some of the  entries of $\hTT$ are equal then they necessarily belong to distinct
columns (in the same tableau), and the corresponding entries of $\TT$ are then chosen
as consecutive integers, increasing from left to right (i.e., with increasing column
indices).
\item[(2)]
Define
$$
\Delta_i\eqdef (\hTT_i - f_{i}\TT - \hTT_{i+1} + f_{i+1}\TT)/2 
\qquad(1\le i\le n), $$
where, by convention, $\hTT_{n+1} \eqdef 1$ and $f_{n+1} \eqdef 0$.
\end{itemize}

\smallskip

\noindent{\bf Example.}
Let $\lala=((3,1),(2,1))$ be a bipartition of $7$, and let
$$
\hTT=\left(
\begin{array}{ccc}
11 & 7 & 3 \\
3  &   &   \\
\end{array}
\;,\;
\begin{array}{cc}
10 & 8 \\
2  &   \\
\end{array}
\right)\in RSSYT\lala.
$$
Computing  $\phi_{\lala}\hTT = (\TT,\Delta)$, the first step yields
$$
\TT=\left(
\begin{array}{ccc}
1 & 4 & 6 \\
5 &   &   \\
\end{array}
\;,\;
\begin{array}{cc}
2 & 3  \\
7 &    \\
\end{array}
\right)\in SYT\lala,
$$
so that $\D\TT=\{1,4,6\}$ and $\N\TT=\{2,3,7\}$.  
Therefore, 
$(f_1\TT,\dots,f_7\TT)=(6,5,5,4,2,2,1)$,
$((\hT^1,\hT^2)_1-f_1\TT,\dots,$ $(\hT^1,\hT^2)_7-f_7\TT)=
(5,5,3,3,1,1,1)$, and
$\Delta=(0,1,0,1,0,0,0)$.

\smallskip

The following claim is easy to verify
\begin{cla}\label{t.phi.B}

\begin{itemize}
\item[{(1)}] $\phi_{(\la^1,\la^2)}$ is a bijection.
\item[{(2)}] If $(\hTT_1,\dots,\hTT_n)$ is the vector of entries of $\hTT$,
in weakly decreasing order, and
$\phi_{(\la^1,\la^2)}(\hT^1,\hT^2) = ((T^1,T^2),\Delta)$, then
$$
\hTT_i -1=f_i\TT+\sum_{j\geq i} 2\Delta_j \qquad(1\le i\le n).
$$
\end{itemize}
\end{cla}
By Claim~\ref{t.phi.B}, for any $\hTT\in RSSYT\lala$
$$
\prod_{i=1}^n q_i^{\hTT_i -1}=\prod_{i=1}^n q_i^{f_i\TT}\cdot
\prod_{j=1}^n (q_1\cdots q_j)^{2\Delta_j}. \leqno(5.2)
$$

Substituting (5.2) into (5.1) we get the formula in the statement of
Lemma~\ref{t.tau.B}.

\qed

\subsection{Proof of Theorem~\ref{t.rep.decomp.B}}

For a permutation $\tau\in B_n$
let the {\it trace} of its action on the polynomial ring $P_n$ be
$$
\tr_{P_n}(\tau) \eqdef
\sum_{m} \langle \tau(m),m\rangle \ \bar q ^{\lambda(m)},
$$
where the sum is over all monomials $m$ in $P_n$,
 $\lambda(m)$
is the exponent partition of the monomial $m$, and the 
inner product is such that the set of all monomials is an orthonormal
basis for $P_n$.

\begin{cla}\label{t.tr.B}
If $\tau\in B_n$ is of cycle type $(\mu^1,\mu^2)$ then the trace of its action on $P_n$ is
$$
\tr_{P_n}(\tau)=\iota[p_{\mu^1,\mu^2}(\bar x,\bar y)],
$$
where $p_{\mu_1,\mu_2}$ is as in Subsection~2.4,
$$
\bar x \eqdef (1,z_1 z_2,z_1 z_2 z_3 z_4,\dots) 
$$
and
$$
\bar y \eqdef (z_1,z_1 z_2 z_3 ,\dots).
$$
\end{cla}
\noindent{\bf Proof.}
Let $\mu^1=(1^{\alpha_1},2^{\alpha_2},\ldots)$ and 
$\mu^2=(1^{\beta_1},2^{\beta_2},\ldots )$.
Then $\langle \tau(m),m\rangle\not=0$ (and is $\pm 1$)
iff all variables $x_i$ with $i$ in the same cycle of $\tau$
have equal exponents in $m$. The sign depends on the number of negative cycles with odd exponent. Thus
$$
\iota^{-1}[\tr_{P_n}(\tau)] =
\prod_{i\ge 1}\left(\sum_{t\ge 0}
(z_1\cdots z_t)^i\right)^{\alpha_i}\cdot
\prod_{j\ge 1}\left(\sum_{s\ge 0}
(-1)^s(z_1\cdots z_s)^j\right)^{\beta_j}=
$$
$$
=p_{\mu^1,\mu^2}(\bar x,\bar y).
$$
\qed

\smallskip

\noindent{\bf Proof of Theorem~\ref{t.rep.decomp.B}.}
From the Straightening Lemma for type $B$ (Corollary \ref{b.2}) it follows that
for any $\tau\in B_n$
$$
\tr_{R_n^B}(\tau)\cdot \prod_{i=1}^n {1\over 1-(q_1\cdots q_i)^2}=
\tr_{P_n}(\tau).
$$
Inserting Claim \ref{t.tr.B}
 and the Frobenius formula for type $B$ 
(see Subsection 2.4)
we obtain
$$
\tr_{R_n^B}(\tau)= \prod_{i=1}^n [1-(q_1\cdots q_i)^2]\ \cdot 
\ \iota[p_{\mu^1,\mu^2}(\bar x,\bar y)]= 
$$
$$
= \prod_{i=1}^n [1-(q_1\cdots q_i)^2]\ \cdot \ \iota[
\sum_{\lambda^1,\lambda^2} \chi^{\lambda^1,\lambda^2}_{\mu^1,\mu^2}
s_{\lambda^1}(\bar x)s_{\lambda^2}(\bar y)].
$$
Applying the linearity of $\iota$ and Lemma \ref{t.tau.B}
 gives
$$
\tr_{R_n^B}(\tau)=\sum_{\la^1,\la^2} \chi^{\la^1,\la^2}_{\mu^1,\mu^2} 
\sum_{\TT\in SYT\lala} \prod_{i=1}^n q_i^{f_i\TT}.
$$

\qed

\section{Combinatorial Identities}

In this section we apply the above algebraic setting
to obtain new combinatorial identities.
The basic algebraic-combinatorial tool here is the Hilbert series
of polynomial rings with respect to multi-degree rearranged into a
weakly decreasing sequence (i.e., a partition).
The computation of the Hilbert series
follows in general the one presented in \cite[\S 8.3]{Re}
for total degree; the major difference is the fact that total degree
gives a {\it grading} of the polynomial ring, whereas rearranged
multi-degree only leads to a {\it filtration}:
$P_n=\cup_\la P_\la^{\underline{\triangleleft}}$ and 
$P_\la^{\underline{\triangleleft}}\cdot P_\mu^{\underline{\triangleleft}}
\subseteq P_{\la+\mu}^{\underline{\triangleleft}}$.
For this reason, not every homogeneous basis for the coinvariant algebra 
will do: the numerators of the RHS in Theorems \ref{7.2} and \ref{7.3} below 
are generating functions for the descent basis (signed descent basis, 
respectively), but other bases (e.g., Schubert polynomials) have different 
generating functions, and are therefore not appropriate for use here. 
The success of the argument relies on properties of the descent basis.

\subsection{Main Combinatorial Results}

For any signed permutation $\sigma\in B_n$ let
$d_{i}(\sigma )$,  
$n_{i}(\sigma )$, and   
$\varepsilon _{i} (\sigma )$ have the same meaning as in 
Subsection 1.2.1. 

The main result of this section is

\begin{thm}\label{7.1}
Let $n \in \bbp$. Then
\[ \sum _{\sigma \in B_{n}} \prod_{i=1}^n 
q_i^{d_i(\sigma)+n_i(\sigma^{-1})} = \sum _{\sigma \in B_{n}}\prod_{i=1}^n 
q_i^{2d_i(\sigma)+\varepsilon_i(\sigma)}.
 \]
\end{thm}

\noindent{\bf Proof.}
Theorem \ref{7.1} is an immediate consequence of Theorems \ref{7.3} 
and \ref{7.4} below.

\qed

\smallskip

For any partition $\lambda=(\la_1,\dots,\la_n)$ with at most $n$ positive parts let 
$$
m_j(\la) \eqdef
|\ \{1\le i\le n\ |\ \la_i=j\}\ |\qquad (\forall j\ge 0),
$$
and let ${n \choose \bar m(\lambda)}$ denote the multinomial coefficient
${n \choose m_0(\lambda),m_1(\lambda), \dots}$.

By considering the Hilbert series of the polynomial ring $P_n$ with respect to weakly decreasing multi-degree we obtain

\begin{thm}\label{7.2} 
Let $n \in \bbp$. 
Then
\[ 
\sum _{\ell(\lambda)\le n} 
{n\choose \bar m(\lambda)} \prod_{i=1}^n q_i^{\lambda_i} = 
\frac{{
\sum _{\pi \in S_{n}} \prod_{i=1}^n q_i^{d_i(\pi)}
}}
{{
 \prod_{i=1}^{n}}(1-q_1\cdots q_i)}
\]
in $\bbz[[q_1, \ldots , q_n]]$, where the sum on the left-hand side is taken over all
partitions with at most $n$ parts. 
\end{thm}

The theorem will be proved using the following lemma.

\smallskip

Recall from Subsection 3.2 the definitions of the index permutation
$\pi(m)$ and the complementary partition $\mu(m)$ of a monomial $m\in P_n$.

\begin{lem}\label{7.2'}
The mapping $m\longmapsto (\pi(m),\mu(m)')$ is a bijection between the set of all monomials
in $P_n$ and the set of all pairs $(\pi,\tilde\mu)$, where
$\pi\in S_n$ and $\tilde\mu$ is a partition with at most $n$ parts.
\end{lem}

\noindent{\bf Proof of lemma \ref{7.2'}.}
The indicated mapping is clearly into the claimed set, since the largest part of the complementary partition $\mu(m)$ is at most $n$.
Conversely, to each pair $(\pi,\tilde\mu)$
as in the statement of the lemma, associate the $\prec$-maximal monomial 
$m$ in the expansion of $a_\pi\cdot e_{\tilde\mu'}$.
Then $\pi(m)=\pi$ and $\mu(m)'=\tilde\mu$. Thus the mapping is a bijection.

\qed

\noindent{\bf Proof of Theorem \ref{7.2}.}
Recall the notation $\la(m)$ for the exponent partition of a monomial $m\in P_n$. 
For any partition with $n$ parts denote by $\bar q^\la$ the product 
$\prod_{i=1}^n q_i^{\la_i}$.

For any partition $\la$ with at most $n$ parts,
${n\choose \bar m(\la)}$ is the number of monomials in $P_n$
with exponent partition equal to $\la$. Therefore
the Hilbert series of the polynomial ring $P_n$ by 
exponent partition
is equal to the LHS of the theorem.

On the other hand, 
since $\la(m)=\la(a_{\pi(m)})+\mu(m)'$, using Lemma \ref{7.2'} we get
$$
\sum_{m\in P_n}\bar q^{\la(m)}=
\sum_{m\in P_n} {\bar q}^{\lambda(a_{\pi(m)})+\mu(m)'}
=\sum_{\pi\in S_n} 
{\bar q}^{\lambda(a_\pi)}\cdot \sum_{\mu}\bar{q}^{\mu},
$$
where $\mu$ in the last sum runs through all partitions having at most 
$n$ parts. By Claim \ref{t.g.des} and Observation~\ref{t.lambda}, 
this product is equal to the RHS of the theorem.

\qed

\smallskip

A well known result, 
attributed by Garsia~\cite{Gar79} to Gessel~\cite{Ge1}, follows. 

\begin{cor}\label{7.g}
Let $n \in \bbp$. Then
\[
\frac{\sum_{\pi \in S_{n}}t^{\des(\pi )
} q^{\maj(\pi )}}
{{
\prod_{i=0}^{n}} (1-tq^{i})} 
=\sum _{r \geq 0} [r+1]_{q}^{n}t^{r}. 
\]
in $\bbz[q][[t]]$.
\end{cor}

\noindent{\bf Proof.} Substitute, in Theorem \ref{7.2}, $q_1=qt$
and $q_2=q_3=\dots=q_n=q$, 
divide both sides by $1-t$,
and apply (3.1) and (3.2). We obtain
$$
\frac{\sum_{\pi \in S_{n}}t^{\des(\pi )
} q^{\maj(\pi )}}
{{
\prod_{i=0}^{n}} (1-tq^{i})}=
\sum_{\ell(\lambda)\le n} {n\choose \bar m(\lambda)} \cdot
q^{\Sigma_i \lambda_i} \cdot {t^{\lambda_1}\over 1-t} =
$$
$$
=\sum_{r=0}^\infty\ \ \sum_{\ell(\lambda)\le n\ ,\  \lambda_1\le r}{n\choose \bar m(\lambda)} q^{\Sigma_i\lambda_i}t^r. 
$$
The coefficient of $t^r$ is
$$
\sum_{\ell(\lambda)\le n\ ,\ \lambda_1\le r} 
{n\choose \bar m(\lambda)} q^{\Sigma_i\lambda_i}=
\sum_{(\ell_1,\dots,\ell_n)\in [0,r]^n} q^{\Sigma_i \ell_i}=
\ (\sum_{j=0}^r q^j)^n\ =\ [r+1]_q^n.
$$

\qed

\smallskip

The same Hilbert series of $P_n$ may be computed in a different way,
by considering the signed descent basis for the coinvariant algebra of type $B$ and applying the Straightening Lemma for this type.

\begin{thm}\label{7.3} 
With notations as in Theorem \ref{7.2}
\[ \sum _{\ell(\lambda)\le n} 
{n\choose \bar m(\lambda)} \prod_{i=1}^n q_i^{\lambda_i} = 
\frac{
{
\sum _{\sigma \in B_{n}} \prod_{i=1}^n 
q_i^{2d_i(\sigma)+\varepsilon_i(\sigma)}
}}
{{
 \prod_{i=1}^{n}}(1-q_1^{2}\cdots q_i^2)}
\]
in $\bbz[[q_1, \ldots , q_n]]$, where the sum on the left-hand side runs through all
partitions with at most $n$ parts.
\end{thm}

\noindent
Recall from Subsection 5.2 the definitions of 
the signed index permutation $\sigma(m)\in B_n$ and
the complementary partition $\mu_B(m)$ of a monomial $m\in P_n$.

\begin{lem}\label{7.3'}
The mapping $m\longmapsto (\sigma(m),\mu_B(m)')$ is a bijection between the set of all monomials
in $P_n$ and the set of all pairs $(\sigma,\tilde\mu)$, where
$\sigma\in B_n$ and $\tilde\mu$ is a partition with at most $n$ parts.
\end{lem}

\noindent{\bf Proof of lemma \ref{7.3'}.}
Similar to the proof of Lemma \ref{7.2'}.
Here to each pair $(\sigma,\tilde\mu)$
associate the $\prec$-maximal monomial 
$m$ in $b_\sigma\cdot e_{\tilde\mu'}(x_1^2,\dots,x_n^2)$,
and then
$\sigma(m)=\sigma$ and $\mu_B(m)'=\tilde\mu$. 

\qed

\noindent{\bf Proof of Theorem \ref{7.3}.}
Again, the LHS of the theorem is equal to 
the Hilbert series of the polynomial ring $P_n$ by exponent partition.
On the other hand, 
since $\la(m)=\la(b_{\sigma(m)})+2\mu_B(m)'$, Lemma \ref{7.3'} gives
$$
\sum_{m\in P_n}\bar q^{\la(m)}=
\sum_{m\in P_n} {\bar q}^{\lambda(b_{\sigma(m)})+2\mu_B(m)'}
=\sum_{\sigma\in B_n} 
{\bar q}^{\lambda(b_\sigma)}\cdot \sum_{\mu}\bar q ^{2\mu},
$$
where $\mu$ in the last sum runs through all partitions
having at most $n$ parts.
By Claim \ref{5.4'}, 
this product is equal to the RHS of the theorem.

\qed

\smallskip

Direct combinatorial arguments imply

\begin{thm}\label{7.4} 
With notations as in Theorem \ref{7.2}
\[ \sum _{\ell( \lambda)\le n} 
{n\choose \bar m(\lambda)} \prod_{i=1}^n q_i^{\lambda_i} = \frac{%
{
\sum _{\sigma \in B_{n}} \prod_{i=1}^n 
q_i^{d_i(\sigma)+n_i(\sigma^{-1})}
}}
{{
 \prod_{i=1}^{n}}(1-q_1^{2}\cdots q_i^2)}
\]
in $\bbz[[q_1, \ldots ,q_n]]$.
\end{thm}

\noindent{\bf Proof.}
Define the subset $T\subseteq B_n$ by
$$
T \eqdef
\{\pi \in B_{n}: \; \des(\pi )=0 \} .
$$

\noindent
It is clear from
our definitions that $d_i( \sigma u) =d_i(u)$ 
and $n_i(u^{-1} \sigma^{-1})=n_i(\sigma^{-1})$
for all $\sigma \in T$, $u \in S_{n}$ and $1\le i\le n$. 
Therefore
\begin{eqnarray*}
\sum _{\pi \in B_{n}} \prod_{i=1}^n 
q_i^{d_i(\pi)+n_i(\pi^{-1})} & = &
\sum _{u \in S_{n}} \sum _{\sigma \in T} 
\prod_{i=1}^n 
q_i^{d_i(\sigma u)+n_i((\sigma u)^{-1})}\\
 & = &
\sum _{u \in S_{n}} \sum _{\sigma \in T} 
\prod_{i=1}^n 
q_i^{d_i(u)+n_i(\sigma^{-1})}\\
 & = &
\sum _{u \in S_{n}}  
\prod_{i=1}^n 
q_i^{d_i(u)}\; \cdot \;
\sum _{\sigma \in T}
\prod_{i=1}^n 
q_i^{n_i(\sigma^{-1})}.\\
\end{eqnarray*}

\noindent
An element $\sigma\in T$ is uniquely determined by the set $\N(\sigma^{-1})$. Hence
$$
\sum_{\sigma \in T} \prod_{i=1}^n q_i^{n_i(\sigma^{-1})}=
 \prod_{i=1}^{n}(1+q_1\cdots q_i) .
$$
Theorem \ref{7.2} completes the proof. 
\qed

\subsection{Negative and Flag Statistics}

In this subsection we present a central result from \cite{ABR} and show that it may
be obtained as a special case of results from the previous subsection.

\smallskip

For any $\sigma\in B_n$ define the
{\it negative descent multiset} by
$$
\DN(\sigma) \eqdef 
\D(\sigma ) \biguplus \N(\sigma^{-1} ) , 
$$
where $\N(\sigma)$ is the set of positions of negative entries in $\sigma$, 
defined in Subsection 1.2.1.
Note that $\DN(\sigma )$ can be defined rather naturally also in purely
Coxeter group theoretic terms. In fact, for $i \in [n]$ let $\eta _{i}
\in B_{n}$ be defined by
\[ 
\eta_{i} \eqdef [1, \ldots , i-1,-i,i+1, \ldots ,n], 
\]
so $\eta_{1} =s_{0}$. Then $\eta_{1}, \ldots, \eta_{n}$ are reflections
 of $B_{n}$
(in the Coxeter group sense; see, e.g., \cite{Hum}).
Clearly
\[
\DN (\sigma ) = 
\{ i \in [n-1]: \; l(\sigma s_{i})< l(\sigma ) \}
\biguplus  
\{ i \in [n]: \; l(\sigma ^{-1}\eta _{i} ) < l(\sigma ^{-1} ) \}  . 
\]
For $\sigma \in B_{n}$ let
\[ \dn (\sigma ) \eqdef |\DN(\sigma )| 
,\qquad
\mn (\sigma ) \eqdef \sum _{i \in \DN(\sigma )} i . \]
Recall the notation
\[ fmaj (\sigma ) =2 \, maj  (\sigma ) + |\N(\sigma )| . \]
and let
$$
 \fd(\sigma ) \eqdef 2 \, \des(\sigma )+\varepsilon_{1}(\sigma ). 
$$
Then the two pairs
 of statistics
$(\fd,\fm)$ and $(\dn,\mn)$ are equidistributed over $B_{n}$.

\smallskip

\begin{cor}\label{t.des.maj.B} \cite[Corollary 4.5]{ABR}
Let $n \in \bbp$. Then
\[ \sum _{\sigma \in B_{n}} t^{\dn(\sigma )}q^{\mn(\sigma )} =
\sum _{\sigma \in B_{n}} t^{\fd(\sigma )} q^{\fm (\sigma )} .  \]
\end{cor}

\noindent{\bf Proof.}
Substitute in Theorem \ref{7.1} : $q_1=qt, q_2=\dots=q_n=q$,
and apply the identities
$$
d_1(\sigma)+n_1(\sigma^{-1})=\dn(\sigma);
$$
$$
\sum_{\i=1}^n (d_i(\sigma)+n_i(\sigma^{-1}))=\mn(\sigma);
$$
and
$$
\sum_{\i=1}^n (2\cdot d_i(\sigma)+\varepsilon_i(\sigma))=\fm(\sigma).
$$

\qed





\smallskip

{\bf Acknowledgments.} 
This paper grew out of stimulating discussions with Dominique Foata and 
Ira Gessel during the conference ``Classical Combinatorics"
in honor of Foata's 65$^{th}$ birthday.
The authors are indebted to Ira Gessel for the idea of using the coinvariant 
algebra in the study of multivariate statistics.
Thanks are also due to Eli Bagno and Richard Stanley for useful comments.


\begin{thebibliography}{xx}

\bibitem{ABR}
R.\ M.\ Adin, F.\ Brenti and Y.\ Roichman,
{\it Descent numbers and major indices for the hyperoctahedral group},
Adv.\ Appl.\ Math.~27 (2001), 210--224.

\bibitem{AR1} 
R.\ M.\ Adin and Y.\ Roichman,
{\it A flag major index for signed permutations},
Proc.\ 11-th Conference on Formal Power Series and Algebraic Combinatorics,
Universitat Polit\`ecnica de Catalunya, Barcelona, 1999, 10--17.

\bibitem{AR}
R.\ M.\ Adin and Y.\ Roichman, 
{\it The flag major index and group actions on polynomial rings}, 
Europ.\ J.\ Combin.~22 (2001), 431--446.

\bibitem{Al1} 
E.\ E.\ Allen, 
{\it The descent monomials and a basis for the diagonally symmetric polynomials}, 
J.\ Alg.\ Combin.~3 (1994), 5--16.


\bibitem{ATY}
S.\ Ariki, T.\ Terasoma and H.\ F.\ Yamada,
{\it Higher Specht polynomials},
Hiroshima Math. J.~27 (1997), 177--188.

\bibitem{Ba} 
H.\ Barcelo, 
{\it Young straightening in a quotient $S_n$-module}, 
J.\ Alg.\ Combin.~2 (1993), 5--23.

\bibitem{BGG} 
I.\ N.\ Bernstein, I.\ M.\ Gelfand and S.\ I.\ Gelfand,
{\it Schubert cells and cohomology of Schubert spaces $G/P$},
Usp.\ Mat.\ Nauk.~28 (1973), 3--26.

\bibitem{BB} 
A.\ Bj\"{o}rner and F.\ Brenti, 
{\it Combinatorics of Coxeter Groups}, 
Graduate Texts in Mathematics, Springer-Verlag, to appear. 

\bibitem{Bre94}
F.\ Brenti, 
{\it $q$-Eulerian polynomials arising from Coxeter groups}, 
Europ.\ J.\ Combin.~15 (1994), 417--441. 




\bibitem{De} 
M.\ Demazure, 
{\it Invariants sym\'etriques entiers des groupes de Weyl et torsion}, 
Invent.\ Math.~21 (1973), 287--301.

\bibitem{Fopc}
D.\ Foata, 
personal communication, July 2000. 

\bibitem{FH}
D.\ Foata and G.\ N.\ Han, 
{\it Calcul basique des permutations signees. I. Longueur et nombre 
d'inversions}, 
Adv.\ Appl.\ Math.~18 (1997), 489--509. 

\bibitem{Gar}
A.\ M.\ Garsia, 
{\it Combinatorial methods in the theory of Cohen-Macaulay rings}, 
Adv.\ Math.~38 (1980), 229--266.

\bibitem{Gar79}
A.\ M.\ Garsia, 
{\it On the ``maj'' and ``inv'' $q$-analogues of Eulerian polynomials}, 
Linear and Multilinear Algebra~8 (1979/80), 21--34.

\bibitem{GP} 
 A.\ M.\ Garsia and C.\ Procesi, 
{\it On certain graded $S_n$-modules and the $q$-Kostka polynomials}, 
Adv.\ Math.~94 (1992), 82--138.

\bibitem{GR} 
A.\ M.\ Garsia and J.\ Remmel,
{\it Shuffles of permutations and the Kronecker product}, 
Graphs and Combinatorics~1 (1985), 217--263.

\bibitem{GS}
A.\ M.\ Garsia and D.\ Stanton,  
{\it Group actions of Stanley-Reisner rings and invariants of permutation groups}, 
Adv.\ Math.~51 (1984), 107--201.

\bibitem{Ge1} 
I.\ M.\ Gessel, 
{\it Generating functions and enumeration of sequences},
Ph.D. Thesis, M.I.T., 1977.

\bibitem{Ge2} 
I.\ M.\ Gessel, 
{\it Multipartite $P$-partitions and inner products of Schur functions}, 
Contemp.\ Math.~34 (1984), 289--302.

\bibitem{Hi}
H.\ L.\ Hiller, 
{\it Geometry of Coxeter Groups}, 
Res.\ Notes in Math.~54, Pitman, Boston, 1982.

\bibitem{Hum} 
J.\ E.\ Humphreys, 
{\it Reflection Groups and Coxeter Groups},
Cambridge Studies in Advanced Mathematics, no.~29,
Cambridge Univ.\ Press, Cambridge, 1990.

\bibitem{KL} 
D.\ Kazhdan and G.\ Lusztig, 
{\it Representations of Coxeter groups and Hecke algebras}, 
Invent.\ Math.~53 (1979), 165--184.

\bibitem{KW} 
W.\ Kraskiewicz and J.\ Weyman, 
{\it Algebra of coinvariants and the action of a Coxeter element},
Bayreuther Math.\ Schriften~63 (2001), 265--284 
(preprint: University of Torun, 1987).

\bibitem{Md} 
I.\ G.\ Macdonald,
{\em Symmetric Functions and Hall Polynomials}, 
second edition, Oxford Math.\ Monographs, 
Oxford Univ.\ Press, Oxford, 1995.

\bibitem{MM} 
P.\ A.\ MacMahon, 
{\em Combinatory Analysis}, 
Chelsea, New York, 1960. 
(Originally published in 2 vols. by Cambridge University Press, 1915-1916.)
 
\bibitem{MY} 
H.\ Morita and H.\ F.\ Yamada, 
{\it Higher Specht polynomials for the complex reflection group $G(r,p,n)$},
Hokkaido Math.\ J.~27 (1998), 505--515.

\bibitem{Rei93a}
V.\ Reiner, 
{\it Signed permutation statistics}, 
Europ.\ J.\ Combin.~14 (1993), 553--567.




\bibitem{Re} 
C.\ Reutenauer, 
{\it Free Lie Algebras},
London Math.\ Soc.\ Monographs, New Series 7, Oxford Univ.\ Press, 1993.

\bibitem{Sa} 
B.\ E.\ Sagan, 
{\it The Symmetric Group: Representations,Combinatorial Algorithms \& 
Symmetric Functions}, 
Wadsworth \& Brooks/Cole, 1991.

\bibitem{So} 
L.\ Solomon,
{\it The orders of the finite Chevalley groups}, 
J.\ Algebra~3 (1966), 376--393.

\bibitem{Sp} 
T.\ A.\ Springer, 
{\it A construction of representations of Weyl groups}, 
Invent.\ Math.~44 (1978), 279--293.

\bibitem{St72}
R.\ P.\ Stanley, 
{\it Ordered Structures and Partitions}, 
Memoirs Amer.\ Math.\ Soc.\, no.~119, 1972.


\bibitem{St79}
R.\ P.\ Stanley, 
{\it Invariants of finite groups and their applications to combinatorics}, 
Bull.\ Amer.\ Math.\ Soc.\ (new series)~1 (1979), 475--511.

\bibitem{St82}
R.\ P.\ Stanley,
{\it Some aspects of group acting on finite posets}, 
J.\ Combin.\ Theory Ser.\ A~32 (1982), 132--161. 

\bibitem{StEC1}
R.\ P.\ Stanley,
{\it Enumerative Combinatorics}, Vol.~1, 
Wadsworth and Brooks/Cole, Monterey, CA, 1986. 

\bibitem{StEC2}
R.\ P.\ Stanley,
{\em Enumerative Combinatorics}, Vol.~2, 
Cambridge Studies in Advanced Mathematics 62,
Cambridge Univ.\ Press, Cambridge, 1999.

\bibitem{Stg} 
R.\ Steinberg, 
{\it On a theorem of Pittie},
Topology~14 (1975), 173--177.


\bibitem{Stem} 
J.\ Stembridge, 
{\it On the eigenvalues of representations of reflection groups and 
wreath products}, 
Pacific J.\ Math.~140 (1989), 353--396.

\bibitem{TY} 
T.\ Terasoma and H.\ F.\ Yamada,
{\it Higher Specht polynomials for the symmetric group},
Proc.\ Japan Acad.\ Ser.~A Math.\ Sci.~69 (1993), 41--44.

\end{thebibliography}
\end{document}